\documentclass[aps,preprint,groupedaddress,showkeys]{revtex4-1}
\usepackage{graphicx}
\usepackage{dcolumn}
\usepackage{bm}
\usepackage{amssymb}
\usepackage{amsthm}
\usepackage{amsmath}
\usepackage{graphicx}
\usepackage{epstopdf}

\newtheorem{theorem}{\hskip\parindent\bf Theorem}
\newtheorem{remark}{\hskip\parindent\bf Remark}
\begin{document}

\preprint{}

\title[Hopf-Hopf bifurcation and chaotic attractors in a delayed diffusive predator-prey model with fear effect]{Hopf-Hopf bifurcation and chaotic attractors in a delayed diffusive predator-prey model with fear effect}

\author{Daifeng Duan \textsuperscript{1}, Ben Niu*\textsuperscript{2}, Junjie Wei \textsuperscript{1}}

\affiliation{\textsuperscript{1}Department of Mathematics, Harbin Institute of Technology, Harbin, Heilongjiang 150001, P.R.China.}
\affiliation{\textsuperscript{2}Department of Mathematics, Harbin Institute of Technology, Weihai 264209, P.R.China.\\*Corresponding author, niubenhit@163.com, niu@hit.edu.cn}

\date{\today}

\begin{abstract}
We investigate a diffusive predator-prey model by incorporating the fear effect into prey population, since the fear of predators could visibly reduce the reproduction of prey.
By introducing the mature delay as bifurcation parameter,
we find this makes the predator-prey system more
complicated and usually induces Hopf and Hopf-Hopf bifurcations.
The formulas determining the properties of Hopf and Hopf-Hopf bifurcations by computing the normal form on the center manifold are given.
Near the Hopf-Hopf bifurcation point we give the
detailed bifurcation set by investigating the universal unfoldings.
Moreover, we show the existence of quasi-periodic orbits on three-torus near a Hopf-Hopf bifurcation point, leading to a strange attractor when further varying the parameter. We also find the existence of Bautin bifurcation numerically, then simulate the coexistence of stable constant stationary solution and periodic solution near this Bautin bifurcation point.
\end{abstract}

                           \maketitle

\section{Introduction}
Studying the predator-prey mechanisms is an important topic on populations, communities and ecosystems. Many mathematical models such as ordinary differential equations, partial differential equations are established to research the growth of population and the spatial distribution law \citep{Creel2008Relationships,Cresswell2011Predation,Su2009Hopf}. In order to rationally exploit biological resources, many scholars have considered
predator-prey models with delay, mainly concentrate on local and global stabilities of equilibria and bifurcations \citep{Faria2001Hopf,Song2005Local,Chakraborty2011bif,Chen2014Bif}.
May \citep{May1973Time} first proposed and discussed the following system with time dealy:
\begin{eqnarray}\label{aa1}
\begin{cases}
\dot{u}(t)=u(t)[r_{1}-au(t-\tau)-p v(t)],\\
\dot{v}(t)=v(t)[-r_{2}+c u(t)-mv(t)],
\end{cases}
\end{eqnarray}
where $u(t)$ and $v(t)$ represent the densities of prey and predator populations at time $t$, respectively.
$\tau$ is the generation time of the prey species, $r_{1}$ denotes the intrinsic growth rate of the prey, $r_{2}$ is the death rate of the predator. $a$ and $m$ represent the death rates due to intra-species competition of the prey and the predator, respectively. $p$ is the capturing rate of the prey by the predator and $c$ is the conversion rate of the prey to the predator. For system (\ref{aa1}),
Song and Wei \citep{Song2005Local} considered the existence of local Hopf bifurcations and paid special attention to the global existence of bifurcating periodic solutions. They showed that the local Hopf bifurcation could be extended to the global Hopf bifurcation under certain conditions.

One of the main trends in the theoretical work of predator-prey dynamics is to derive more realistic models to explain complicated biological phenomena.
A recent research on animals on land showed that the cost of fear changed anti-predator defences, thus it greatly reduced the reproduction of prey. The previous view is that predators can only influence prey populations by killing directly. In the wild, it is easy to observe that the reduction of prey is due to the direct killing of predators \citep{Lima1998Nonlethal}.
Even if the functional responses such as Holling type I-III \citep{Holling1965,Kooij1997qualitative,Yi2009Bif,Liu2018per}, Beddington-DeAngelis \citep{Beddington1975Mutual,Deangelis1975A,Hwang2003Global} are added to the predator-prey models, they can only reflect direct killing.
A new study, however, suggests that the presence of a predator can also change prey behavior and have a greater impact than direct killing \citep{Lima2010Predators,Zanette2011per}.
In addition, the fear of predators may affect the physiological status of infant prey individuals, and is detrimental to their survival.

Recently Wang et.al. \citep{Wang2016Modelling} proposed and analyzed
a two-dimensional autonomous differential equation modeling the fear effect in predator-prey interactions, biologically introduced in
\citep{Creel2008Relationships,Cresswell2011Predation,Lima2010Predators,Zanette2011per}. They considered the following model incorporating the cost of fear (indirect effects).
\begin{eqnarray}\label{105a}
\begin{cases}
\dot{u}=u\big[r_{0}f(K,v)-d-a u \big]-g(u) v,\\
\dot{v}=v(-r_{2}+c g(u)),
\end{cases}
\end{eqnarray}
where $r_{0}$ is the birth rate of the prey, $d$ is the natural death rate of the prey. They multiplied the production
term by a fear factor $f(K,v)$ in the first equation of system (\ref{aa1}), where the parameter $K$ reflects the level of fear that prey drives predator. Here $g(u)$ is the functional response of predators to prey densities \citep{Wang2016Modelling}.
Wang et al. \citep{Wang2016Modelling} have studied system (\ref{105a}) with either the linear ($g(u)=pu$) or the Holling type II ($g(u)=p/(1+qu)$) functional response. Their theoretical results suggest that fear effect could stabilize the predator-prey system
by eliminating the existence of periodic solutions. Compared with the classic predator-prey model where Hopf bifurcations are usually supercritical, Hopf bifurcations in their model could be supercritical or subcritical in a series of numerical simulations.

In fact, there usually exist unevenly distributed predators and preys in different spatial locations, which will move or spread to areas with lower population density or abundant food in order to obtain fine environments. Mathematical analysis shows that the predator-prey system with diffusion will exhibit complex characteristics. The dynamical properties mainly include that diffusion coefficients could lead to spatially non-homogeneous bifurcating periodic solutions or Turing instability \citep{Shen2016,Yang2016,Zhao2014,Zheng2014}.

Motivated by these pioneer works, we choose the particular form for the fear effect term $f(K,v)=1/(1+Kv)$ and linear functional response $g(u)=pu$, then add diffusion terms as
\begin{eqnarray}\label{105b}
\begin{cases}
\frac{\partial }{\partial t}u(x,t)=d_{1} \Delta u(x,t)+u(x,t)\Big[\frac{r_{0}}{1+K v(x,t)}-d-au(x,t-\tau)-p v(x,t)\Big]
, \\ 
\frac{\partial}{\partial t} v(x,t)=d_{2} \Delta v(x,t)+v(x,t)\big[-r_{2}+cu(x,t)-m v(x,t)\big],
~x\in\Omega, t>0, \\
\frac{\partial u(x,t)}{\partial \overrightarrow{n}}
=\frac{\partial v(x,t)}{\partial \overrightarrow{n}}=0, ~~x\in\partial\Omega, t>0,\\
u(x,t)=u_{0}(x,t)\geq0, v(x,t)=v_{0}(x,t)\geq0, ~x\in\overline{\Omega},t\in[-\tau,0],
\end{cases}
\end{eqnarray}
where $u(x,t)$ and $v(x,t)$ are prey and predator densities respectively,
at location $x\in\Omega$ and time $t>0$, $\Omega$ is a bounded domain in
$\mathbb{R}^{n}$ with a smooth boundary $\partial\Omega$ for $n\geq1$.
$\tau$ represents the generation time of $u(x,t)$,
$d_{1}, d_{2}>0$ denote the diffusion coefficients of prey and predator, respectively, and $\overrightarrow{n}$ is the outward unit normal vector at $\partial\Omega$. The parameters $r_{0}, r_{2}, K, d, a, p, c$ and $m$ are all positive. A homogeneous Neumann boundary condition is imposed so that the population
movement across the boundary is algebraically zero.

We find the delay makes the predator-prey system (\ref{105b}) more complicated and usually induces
stability switches or Hopf-Hopf bifurcations.
With the extended Hopf and Hopf-Hopf bifurcations results and normal form methods \citep{Du2018Double,Faria2000Normal,Faria1995Normal}, we investigate the effect of birth rate of prey and time delay in view of the existence of Hopf-Hopf bifurcation. The main achievements of this paper is a detailed bifurcation analysis about the positive constant stationary solution
of (\ref{105b}) with the one-dimensional spatial domain.
The existence of periodic solutions and quasi-periodic solutions are obtained by using the normal form theory. Moreover, we obtain the existence of strange attractor in term of \citep{Newhouse1993Occurrence}.
Guided by these theoretical analyses, the above results are illustrated numerically.

The rest of this paper is organized as follows.
In Sec. \ref{105d}, we analyze the critical conditions of Hopf and Hopf-Hopf bifurcations.
In Sec. \ref{a1019}, we study the direction and stability of Hopf bifurcation near the positive equilibrium, then derive the third-order truncated normal form for the Hopf-Hopf bifurcation by \citep{Du2018Double}.
In Sec. \ref{a1112}, we show the bifurcation set and phase portraits of system (\ref{105b}).
Besides, some other numerical examples are performed, such as the coexistence of stable constant stationary solution and spatially homogeneous stable periodic solution (see \citep{Guo2018}), the existence of a 2-torus or 3-torus and the transition from quasi-periodic oscillations to chaos.
Some direct calculations are relegated to the appendix A and B.
Throughout this paper, we denote by $\mathbb{N}$ the set of all positive integers, and
$\mathbb{N}_{0} = \mathbb{N} \cup \{0\}$.
\section{Stability of the positive equilibrium and bifurcations in the delayed model}\label{105d}
In the rest part of this paper, we consider the
predator-prey dynamics in one dimensional space, that is $\Omega=(0,l\pi)$, $l\in\mathbb{R}^{+}$.
In terms of a real ecosystem, it can represent longitude or latitude, marine turbulence, the depth of water, etc.

For convenience of the readers we copy (\ref{105b}) here as,
\begin{eqnarray}\label{a}
\begin{cases}
\frac{\partial }{\partial t}u(x,t)=d_{1} \Delta u(x,t)+u(x,t)\Big(\frac{r_{0}}{1+K v(x,t)}-d-au(x,t-\tau)-p v(x,t)\Big), \\
\frac{\partial}{\partial t} v(x,t)=d_{2} \Delta v(x,t)+v(x,t)\big(-r_{2}+cu(x,t)-m v(x,t)\big), ~x\in\Omega, t>0,\\
\frac{\partial u(x,t)}{\partial \overrightarrow{n}}
=\frac{\partial v(x,t)}{\partial \overrightarrow{n}}=0, ~~x=\partial\Omega, t>0.
\end{cases}
\end{eqnarray}
We can see that system (\ref{a}) has an extinct equilibrium $(0,0)$ and a boundary equilibrium $((r_{0}-d)/a,0)$ under the condition $r_{0}>d$.
Besides, when
\begin{eqnarray*}
\mathrm{(H0)}~~~ar_{2}+cd-c r_{0}<0
\end{eqnarray*}
holds, (\ref{a}) has a positive constant (spatially homogeneous) stationary solution $E_{\ast}=(u_{\ast},v_{\ast})$, with
\begin{eqnarray}\label{b}
\begin{split}
v_{\ast}&=\frac{1}{2K(am+cp)}\Big[-(am+cp+K(ar_{2}+cd))\\
&+\sqrt{(am+cp+K(ar_{2}+cd))^{2}-4K(am+cp)(ar_{2}+cd-c r_{0})}\Big],\\
u_{\ast}&=(r_{2}+m v_{\ast})/c.
\end{split}
\end{eqnarray}
Linearizing (\ref{a}) around $E_{\ast}$ takes the
form
\begin{eqnarray}\label{c}
\frac{\partial}{\partial t}\left(
  \begin{array}{c}
     u(x,t) \\
     v(x,t) \\
  \end{array}
\right)=D  \left(
           \begin{array}{c}
           \Delta u(x,t) \\
            \Delta v(x,t) \\
           \end{array}
         \right)+\tilde{A}\left(
                    \begin{array}{c}
                      u(x,t) \\
                      v(x,t) \\
                    \end{array}
                  \right) +\tilde{B}\left(
                    \begin{array}{c}
                      u(x,t-\tau) \\
                      v(x,t-\tau) \\
                    \end{array}\right),
\end{eqnarray}
where
\begin{equation*}
   D=\left(\begin{array}{cc}
d_{1}  & 0  \\
0 & d_{2}   \\
\end{array}\right),
   \tilde{A}=\left(\begin{array}{cc}
0  & a_{12}  \\
a_{21} & a_{22}   \\
\end{array}\right),
  \tilde{B}=\left(\begin{array}{cc}
b_{11}  & 0  \\
0 & 0   \\
\end{array}\right),
\end{equation*}
and
\begin{eqnarray}\label{b629}
\begin{split}
a_{12}=\frac{-K r_{0}u_{\ast}}{(1+K v_{\ast})^{2}}-pu_{\ast}<0,~a_{21}=cv_{\ast}>0,
~a_{22}=-m v_{\ast}<0,~b_{11}=-au_{\ast}<0.
\end{split}
\end{eqnarray}

As the Laplace operator has eigenvalues $-k^{2}/l^{2}$~
$(k = 0, 1, ... )$, with corresponding normalized eigenfunctions
\begin{eqnarray}\label{a706}
\begin{split}
\gamma_{k}(x) = \frac{\mathrm{cos}\frac{k}{l}x}{\|\mathrm{cos}\frac{k}{l}x\|_{L^{2}}}=\left\{
\begin{array}{l}
\sqrt{\frac{1}{l\pi}}, ~~~~~~~~~k=0,\\
\sqrt{\frac{2}{l\pi}}\mathrm{cos}\frac{k}{l}x,~~k\geq1,
\end{array}
\right.
\end{split}
\end{eqnarray}
$\lambda$ is a characteristic value of (\ref{c}) if and only if for some $k = 0,1, ... $, $\lambda$ satisfies
\begin{equation*}
   \mathrm{det}\left(\begin{array}{cc}
\lambda+d_{1}\frac{k^{2}}{l^{2}}-b_{11}e^{-\lambda \tau}  & -a_{12}  \\
-a_{21} &  \lambda+d_{2}\frac{k^{2}}{l^{2}}-a_{22}  \\
\end{array}\right)=0.
\end{equation*}
That is, $\lambda$ solves
\begin{eqnarray}\label{e1}
  \lambda^{2}+A_{k}\lambda+B_{k}+(-b_{11}\lambda+C_{k})e^{-\lambda\tau}=0, ~k\geq0,
\end{eqnarray}
with
$A_{k}=(d_{1}+d_{2})\frac{k^{2}}{l^{2}}-a_{22}>0$,
$B_{k}=d_{1}d_{2}\frac{k^{4}}{l^{4}}-a_{22}d_{1}\frac{k^{2}}{l^{2}}-a_{12}a_{21}>0$,
$C_{k}=-b_{11}d_{2}\frac{k^{2}}{l^{2}}+a_{22}b_{11}>0$.
We find that $\lambda=0$ is not the root of (\ref{e1}) since $B_{k}+C_{k}>0$, it is impossible for system (\ref{a}) to possess Turing instability.
When $\tau=0$, the characteristic equation (\ref{e1}) becomes
$\lambda^{2}+(A_{k}-b_{11})\lambda+B_{k}+C_{k}=0$.
Obviously, $A_{k}-b_{11}>0$ and $B_{k}+C_{k}>0$ result in the following conclusion.
\begin{theorem}\label{105c}
When $\tau=0$, the positive equilibrium $E_{\ast}$ of system (\ref{a}) is locally asymptotically stable.
\end{theorem}

Applying the same analytical steps as these in Ruan and Wei \citep{Ruan2003On},
when $\tau>0$, letting $\lambda=i\omega$~$(\omega>0)$ in (\ref{e1}), we obtain
\begin{eqnarray}\label{623a}
\left\{
\begin{array}{l}
 \mathrm{cos}\omega\tau =\frac{d_{1}b_{11}\omega^{2}k^{2}/l^{2}-B_{k}C_{k}}{b_{11}^{2}\omega^{2}+C_{k}^{2}}=C_{k}~(\omega) (<0), \\
  \mathrm{sin}\omega\tau = \frac{A_{k}C_{k}\omega+B_{k}b_{11}\omega-b_{11}\omega^{3}}{b_{11}^{2}\omega^{2}+C_{k}^{2}}=S_{k}(\omega),
\end{array}
\right.
\end{eqnarray}
from which it follows that
\begin{eqnarray}\label{f1}
\omega^{4}+(A_{k}^{2}-2B_{k}-b_{11}^{2})\omega^{2}+B_{k}^{2}-C_{k}^{2}=0,
\end{eqnarray}
the roots of (\ref{f1}) are
\begin{eqnarray}\label{c707}
 \omega_{k}^{\pm} =\bigg[\frac{1}{2}\Big(-A_{k}^{2}+2B_{k}+b_{11}^{2}\pm \sqrt{(A_{k}^{2}-2B_{k}-b_{11}^{2})^{2}-4(B_{k}^{2}-C_{k}^{2})}\Big)\bigg]^{1/2}.
\end{eqnarray}
We make the following assumptions,
\begin{eqnarray*}
&&\mathrm{(H1)}~ \{-A_{k}^{2}+2B_{k}+b_{11}^{2}<0 ~\mathrm{and} ~B_{k}-C_{k}>0\}, \mathrm{or}~ (A_{k}^{2}-2B_{k}-b_{11}^{2})^{2}<4(B_{k}^{2}-C_{k}^{2}).\\
&&\mathrm{(H2)}~ B_{k}-C_{k}<0.\\
&&\mathrm{(H3)}~ B_{k}-C_{k}>0, -A_{k}^{2}+2B_{k}+b_{11}^{2}>0 ~\mathrm{and}~ (A_{k}^{2}-2B_{k}-b_{11}^{2})^{2}>4(B_{k}^{2}-C_{k}^{2}).
\end{eqnarray*}
If (H1) is satisfied, Eq.(\ref{f1}) has no positive root, then Eq.(\ref{e1}) has no purely imaginary root.
If (H2) is satisfied, Eq.(\ref{f1}) has one positive root $\omega_{k}^{+}$, then Eq.(\ref{e1}) has a pair of purely imaginary roots $\pm i\omega_{k}^{+}$ at $\tau_{k}^{j+}$, with
\begin{eqnarray*}
\tau_{k}^{j+}=
 \frac{\pi-\mathrm{arcsin} S_{k}(\omega_{k}^{+})+2j\pi}{\omega_{k}^{+}},~j\in\mathbb{N}_{0}.
\end{eqnarray*}
If (H3) is satisfied, Eq.(\ref{f1}) has two positive roots $\omega_{k}^{\pm}$, then Eq.(\ref{e1}) has two pairs of purely imaginary roots $\pm i\omega_{k}^{\pm}$ at $\tau_{k}^{j\pm}$, $j\in \mathbb{N}_{0}$, with
\begin{eqnarray}\label{603h}
\tau_{k}^{j\pm}=
 \frac{\pi-\mathrm{arcsin} S_{k}(\omega_{k}^{\pm})+2j\pi}{\omega_{k}^{\pm}},~j\in\mathbb{N}_{0}.
\end{eqnarray}
If (H0) and (H3) hold for some $k\in \mathbb{N}_{0}$, differentiating the two sides of (\ref{e1}) with respect to $\tau$, after some basic calculations similar with those in \citep{Zhao2015Dynamics},
we have
\begin{eqnarray*}
\mathrm{Re}\bigg(\frac{d\lambda}{d\tau}\bigg)\bigg|_{\tau=\tau_{k}^{j\pm}}=\pm \frac{\sqrt{(A_{k}^{2}-2B_{k}-b_{11}^{2})^{2}-4(B_{k}^{2}-C_{k}^{2})}}{C_{k}^{2}+b_{11}^{2}\omega_{k}^{\pm2}}.
\end{eqnarray*}
Thus $\mathrm{Re}\big( \frac{d\lambda}{d\tau}\big)\big|_{\tau=\tau_{k}^{j+}}>0$, $\mathrm{Re}\big( \frac{d\lambda}{d\tau}\big)\big|_{\tau=\tau_{k}^{j-}}<0$ for $j\in \mathbb{N}_{0}$.

Denote $\mathcal{D}_{1}=\{k\in \mathbb{N}_{0}: \mathrm{(H2)~holds}\}$, $\mathcal{D}_{2}=\{k\in \mathbb{N}_{0}: \mathrm{(H3)~holds}\}$.
Note that $\mathcal{D}_{1}$ is a finite set, since $\lim\limits_{k\rightarrow\infty}(B_{k}-C_{k})=+\infty$. $\mathcal{D}_{2}$ is also a finite set, since $\lim\limits_{k\rightarrow\infty}(-A_{k}^{2}+2B_{k}+b_{11}^{2})=-\infty$.
Obviously, $\tau_{k}^{j\pm}$ is monotonically
increasing on $j$ for fixed $k\in \mathcal{D}_{2}$, so $\tau_{k}^{0\pm}=\mathrm{min}_{j\in \mathbb{N}_{0}}\{\tau_{k}^{j\pm}\}$ for fixed $k$.
We define the smallest $\tau$ such that the stability of $E_{\ast}$ may change:
\begin{eqnarray*}
\bar{\tau}\overset{\mathrm{def}}{=}\left\{
\begin{array}{ll}
\mathrm{min}\{\tau_{k}^{0+}, ~k\in \mathcal{D}_{1}\}, & \mathrm{if}~\mathrm{(H2)}~\mathrm{holds},\\
\mathrm{min}\{\tau_{k}^{0+}, \tau_{k}^{0-}, ~k\in\mathcal{D}_{2}\}, & \mathrm{if}~\mathrm{(H3)}~\mathrm{holds} . \nonumber
\end{array}
\right.
\end{eqnarray*}

From the discussion above we know that, under the assumption (H3), there are two sequences value of $\tau$, $\{\tau_{k}^{j+}\}$ and $\{\tau_{k}^{j-}\}$, such that (\ref{e1}) has a pair
of purely imaginary roots when $\tau=\tau_{k}^{j\pm}$, respectively.
Assume that $\tau_{1}, \tau_{2}, \cdots, \tau_{m}$ ($m\in\mathbb{N}$) are chosen by
$\{\tau_{k}^{j\pm}, k\in\mathcal{D}_{2},j\in\mathbb{N}_{0}\}$ and $\bar{\tau}<\tau_{1}< \cdots< \tau_{m}$.
We have the following results.
\begin{theorem}\label{a912}
Suppose (H0) holds, for system (\ref{a}), the following statements hold true.\\
$\mathrm{(i)}$~If (H1) holds for any $k\in \mathbb{N}_{0}$, then the positive equilibrium $E_{\ast}$ is locally asymptotically stable for all $\tau\geq0$.\\
$\mathrm{(ii)}$~If (H2) holds, then
the positive equilibrium $E_{\ast}$ is locally asymptotically stable for $\tau\in[0,\bar{\tau})$,
and unstable for $\tau>\bar{\tau}$. Besides, system (\ref{a}) undergoes a Hopf bifurcation at
the positive equilibrium $E_{\ast}$ when $\tau=\tau_{k}^{j+}$ for $j\in \mathbb{N}_{0}$, $k\in\mathcal{D}_{1}$.\\
$\mathrm{(iii)}$~If (H3) holds, there exist finite critical points $\bar{\tau}<\tau_{1}< \cdots< \tau_{m}$, and $\tau_{1}, \tau_{2},\cdots,\tau_{m}\in\{\tau_{k}^{j\pm}, k\in\mathcal{D}_{2},j\in\mathbb{N}_{0}\}$ for $m\in\mathbb{N}_{0}$ such that when $\tau\in[0,\bar{\tau})\cup(\tau_{1},\tau_{2})\cup\cdots\cup(\tau_{m-1},\tau_{m})$,
$E_{\ast}$ is locally asymptotically stable; when
$\tau\in(\bar{\tau},\tau_{1})\cup(\tau_{2},\tau_{3})\cup\cdots\cup(\tau_{m},+\infty)$,
$E_{\ast}$ is unstable.
Besides, system (\ref{a}) undergoes a Hopf bifurcation at
the positive equilibrium $E_{\ast}$ when $\tau=\tau_{k}^{j\pm}$ for $j\in \mathbb{N}_{0}$, $k\in\mathcal{D}_{2}$.
\end{theorem}

In this paper, we choose the birth rate $r_{0}$ and time delay $\tau$ as Hopf-Hopf  bifurcation parameters and obtain the following result:
\begin{theorem}\label{b707}
Under assumptions $\mathrm{(H0)}$ and $\mathrm{(H3)}$, if there exists $r_{0}=r_{0}^{\ast}$ such that $\tau_{k_{1}}^{j_{1}+}=\tau_{k_{2}}^{j_{2}-}=\tau^{\ast}$ for some $j_{1}, j_{2}, k_{1}, k_{2}\in\mathbb{N}_{0}$, then Eq.(\ref{e1}) has two pairs of purely imaginary roots $\pm i\omega_{k}^{\pm}$ ($k\in\mathcal{D}_{2}$), when $(\tau,r_{0})=(\tau^{\ast},r_{0}^{\ast})$.
\end{theorem}

\begin{remark}\label{d707}
Under Theorem \ref{b707},
system (\ref{a}) may undergo a Hopf-Hopf bifurcation at the positive constant stationary solution $E_{\ast}$ when $(\tau,r_{0})=(\tau^{\ast},r_{0}^{\ast})$. In order to better analyze
the dynamical behavior near the bifurcation point, we will give the bifurcation analysis of $E_{\ast}$ about system (\ref{a}).
\end{remark}

\section{Bifurcation analysis}\label{a1019}
In this section, we will investigate the direction and stability of Hopf bifurcation
near the positive equilibrium $E_{\ast}$ by applying the center manifold theorem and normal form
theory of the partial differential equations presented in Faria \citep{Faria2000Normal} and Wu \citep{Wu1996Theory}. Besides, we will derive the third-order truncated normal form
for the  Hopf-Hopf bifurcation by applying the normal form theory in \citep{Du2018Double}.
\subsection{Hopf bifurcation}
In this subsection, we will obtain the stability of the bifurcating periodic
solutions and show more detailed information of Hopf bifurcation by
using the normal form theory due to \citep{Faria2000Normal}.
Let $\hat{u}(x,t)=u(x,\tau t)-u_{\ast}$, $\hat{v}(x,t)=v(x,\tau t)-v_{\ast}$, and drop the hats for simplification of notation, then system (\ref{a}) can be rewritten as
\begin{eqnarray}\label{a916}
\begin{cases}
\frac{\partial u(x,t)}{\partial t}=
\tau[d_{1}\mathrm{\Delta} u(x,t)+a_{12}v(x,t)+b_{11}u(x,t-1)+\alpha_{1}u(x,t)v(x,t)\\
~~-au(x,t)u(x,t-1)+\alpha_{2}v^2(x,t)+\alpha_{3}v^3(x,t)+\alpha_{4}u(x,t)v^2(x,t)]+\mathrm{h.o.t.},\\
\frac{\partial v(x,t)}{\partial t}=\tau[d_{2}\mathrm{\Delta } v(x,t)+a_{21}u(x,t)+a_{22}v(x,t)+cu(x,t)v(x,t)-mv^{2}(x,t)],
\end{cases}
\end{eqnarray}
where $a_{12}, a_{21}, a_{22}, b_{11}$ are given by (\ref{b629}), and
\begin{eqnarray}
\begin{split}
\alpha_{1}=\frac{-Kr_{0}}{(1+Kv_{\ast})^{2}}-p,~\alpha_{2}=\frac{K^{2}r_{0}u_{\ast}}{(1+Kv_{\ast})^{3}},~ \alpha_{3}=\frac{-K^{3}r_{0}u_{\ast}}{(1+Kv_{\ast})^{4}},~
\alpha_{4}=\frac{K^{2}r_{0}}{(1+Kv_{\ast})^{3}}.
\end{split}
\end{eqnarray}
Define the real-valued Hilbert space
\begin{eqnarray*}
  X:=\bigg\{(u,v)\in H^{2}(0,l\pi)\times H^{2}(0,l\pi):\bigg(\frac{\partial u}{\partial x},\frac{\partial v}{\partial x}\bigg)\bigg|_{x=0,l\pi}=0\bigg\},
\end{eqnarray*}
and the corresponding complexification $X_{\mathbb{C}}$ has the form
$X_{\mathbb{C}}:=\{x_{1}+ix_{2}, x_{1},x_{2}\in X \}$,
and the complex-valued $L^{2}$ inner product is given by
$\langle \tilde{a}, \tilde{b} \rangle=\int_{0}^{l\pi}(\bar{a}_{1}b_{1}+\bar{a}_{2}b_{2})dx$, for $\tilde{a}=(a_{1},a_{2})^{T}$, $\tilde{b}=(b_{1},b_{2})^{T}\in X_{\mathbb{C}}$.
Define the phase space with the sup norm $\mathcal{C}:=C([-1,0],X_{\mathbb{C}})$, and  write $\varphi_{t}\in \mathcal{C}$ ,
$\varphi_{t}(\theta)=\varphi(t+\theta)$ for $-1\leq\theta\leq0$.
Let $\beta_{k}^{(1)}(x)=(\gamma_{k}(x),0)^{T}$, $\beta_{k}^{(2)}(x)=(0,\gamma_{k}(x))^{T}$ and $\beta_{k}=\{\beta_{k}^{(1)}(x), \beta_{k}^{(2)}(x)\}$, where
$\{\beta_{k}^{(i)}(x)\} (i=1,2)$ is an
orthonormal basis of $X$. Define the subspace of $\mathcal{C}$, that is,
$\mathcal{B}_{k}:=\mathrm{span}\{\langle \phi(\cdot),\beta_{k}^{(j)}\rangle\beta_{k}^{(j)}|\phi\in \mathcal{C},j=1,2\},~k=0,1,2,\cdots$.
Choose
\[  \eta^{k}(\theta,\bar{\tau}) =\begin{cases}
\bar{\tau}(-\frac{k^{2}}{l^{2}}D+\tilde{A})& \theta=0,\\
0 & \theta\in(-1,0),\\
-\bar{\tau} \tilde{B} &  \theta=-1.
\end{cases} \]
By the Riesz representation theorem, there exists a matrix function $\eta^{k}(\theta,\bar{\tau})$ of the bounded variation
for $\theta\in[-1,0]$, such that $-\bar{\tau}D\frac{k^{2}}{l^{2}}\phi(0)+\bar{\tau}L(\phi)=\int_{-1}^{0}d\eta^{k}(\theta,\bar{\tau})\phi(\theta)$, $\phi\in\mathcal{C}$.
The bilinear form on $\mathcal{C}^\ast\times\mathcal{C}$ is defined by
\begin{eqnarray}\label{b916}
(\psi,\phi)=\psi(0)\phi(0)-\int^{0}_{-1}\int^{\theta}_{\xi=0}
\psi(\xi-\theta)\mathrm{d}\eta^{k}(\theta,\bar{\tau})\phi(\xi)\mathrm{d}\xi,
\end{eqnarray}
where $\phi\in\mathcal{C}$, $\psi\in\mathcal{C}^\ast$.
Let $\tau=\bar{\tau}+\mu$, from the discussion in section \ref{105d}, we know that when $\mu=0$ system (\ref{a916}) undergoes a Hopf bifurcation at the equilibrium $(0,0)$, Eq.(\ref{e1}) has a pair of purely imaginary roots $\pm i\omega_{n_{0}}$.
Let $\mathcal{A}$ denote the infinitesimal generator of the semigroup with $\mu=0$ and $n=n_{0}$, and $\mathcal{A}^{\ast}$ denote the formal adjoint of $\mathcal{A}$ under the bilinear form (\ref{b916}).
Let $p(\theta)=p(0)e^{i\omega_{n_{0}}\bar{\tau}\theta}(\theta\in[-1,0])$,
$q(s)= q(0)e^{-i\omega_{n_{0}}\bar{\tau}s}$ $(s\in[0,1])$
be the eigenvectors of $\mathcal{A}$ and $\mathcal{A}^{\ast}$  corresponding to the eigenvalue $i\omega_{n_{0}}\bar{\tau}$.
By direct calculations, we choose
$p(0)=(1,p_{1})^T$, $q(0) = M(1,q_{2})$,
where
$p_{1}=a_{21}/(i\omega_{n_{0}}+d_{2}n_{0}^{2}/l^{2}-a_{22})$,
$q_{2}=a_{12}/(i\omega_{n_{0}}+d_{2}n_{0}^{2}/l^{2}-a_{22})$,
$M=(1+p_{1}q_{2}+\bar{\tau}b_{11}e^{-i\omega_{n_{0}}\bar{\tau}})^{-1}$.

According to the general theory in  \citep{Hassard1981Theory},
we compute the normal form up to the third order for system (\ref{a916}). We leave the detailed procedure in appendix A, where two key values $\mu_{2}$ and $\beta_{2}$ are calculated with
\begin{eqnarray*}
\mu_{2}=-\frac{\mathrm{Re}(c_{1}(0))}{\mathrm{Re}(\lambda'(\bar{\tau}))},
~\beta_{2}=2\mathrm{Re}(c_{1}(0)),
\end{eqnarray*}
and
\begin{eqnarray}\label{c1021}
c_{1}(0)=\frac{i}{2\omega_{n_{0}}\bar{\tau}}\Big
(g_{11}g_{20}-2|g_{11}|^{2}-\frac{|g_{02}|^{2}}{3}\Big)+\frac{g_{21}}{2}.
\end{eqnarray}
We have the following conclusion.
\begin{theorem}\label{a1113}
For system (\ref{a}), the Hopf bifurcation at $\tau=\bar{\tau}$ is supercritical (or subcritical) if $\mu_{2}>0$ (or $\mu_{2}<0$), and the bifurcating
solution is orbitally asymptotically stable (or unstable).
\end{theorem}

\subsection{Hopf-Hopf bifurcation}
Combining theorem \ref{b707} with remark \ref{d707} in section \ref{105d}, Hopf-Hopf bifurcation occurs with the existence of stability switches. In this subsection, we give the universal unfoldings near the Hopf-Hopf bifurcation by applying the normal form theory in \citep{Du2018Double}.

In order to study the qualitative behavior near the critical point $(\tau^{\ast},r_{0}^{\ast})$,
set $\tau=\tau^{\ast}+\mu_{1}$, $r_{0}=r_{0}^{\ast}+\mu_{2}$,
$u(t)=u(\cdot,t)$, $v(t)=v(\cdot,t)$
and $U(t)=(u(t),v(t))^{T}$, $U_{t}(\cdot)=U(t+\cdot)$,
then system (\ref{a916}) can be written as an abstract differential equation
in the phase space $\mathcal{C}:=C([-1,0],X_{\mathbb{C}})$ as follows:
\begin{eqnarray}\label{a1}
\frac{d }{d t}U(t)&=\tau^{\ast}D\Delta U(t)+L_{0}U_{t}+\tilde{F}(\mu_{1},\mu_{2},U_{t}),
\end{eqnarray}
where
$L_{0}U_{t}=\tau^{\ast}(\tilde{A}_{0}U_{t}(0)+\tilde{B}_{0}U_{t}(-1))$,
$\tilde{F}(\mu_{1},\mu_{2},U_{t})=\mu_{1}D\Delta U(t)+
\mu_{1}L_{11}U_{t}+\mu_{2}L_{12}U_{t}+F(0,0,U_{t})$,
$L_{11}U_{t}=\tilde{A}_{0}U_{t}(0)+\tilde{B}_{0}U_{t}(-1)$,
$L_{12}U_{t}=\tau^{\ast}(\tilde{A}_{1}U_{t}(0)+\tilde{B}_{1}U_{t}(-1))$.

When $r_{0}=r_{0}^{\ast}$, differing from (\ref{b}),
denote positive constant stationary solution $(\tilde{u}_{\ast},\tilde{v}_{\ast})$ by
\begin{eqnarray*}
\begin{split}
\tilde{v}_{\ast}&=\frac{1}{2K(am+cp)}\Big[-(am+cp+K(ar_{2}+cd))\\
&+\sqrt{(am+cp+K(ar_{2}+cd))^{2}-4K(am+cp)(ar_{2}+cd-c r^{\ast}_{0})}\Big],\\
\tilde{u}_{\ast}&=(r_{2}+m \tilde{v}_{\ast})/c.
\end{split}
\end{eqnarray*}
and
\begin{eqnarray}\label{610a}
   \tilde{A}_{0}=\left(\begin{array}{cc}
0  & a_{12}^{\ast}  \\
a_{21}^{\ast} & a_{22}^{\ast}   \\
\end{array}\right),
   \tilde{B}_{0}=\left(\begin{array}{cc}
b_{11}^{\ast}  & 0  \\
0 & 0   \\
\end{array}\right),
\end{eqnarray}
with $a_{12}^{\ast}=-\frac{K r_{0}^{\ast}\tilde{u}_{\ast}}{(1+K \tilde{v}_{\ast})^{2}}-p\tilde{u}_{\ast}$,
 $a_{21}^{\ast}=c\tilde{v}_{\ast}$,
$a_{22}^{\ast}=-m \tilde{v}_{\ast}$, $b_{11}^{\ast}=-a\tilde{u}_{\ast}$.
Let $(u_{\ast}',v_{\ast}')$ denote the derivative of $(u_{\ast},v_{\ast})$ with respect to $r_{0}$ evaluated at $r_{0}^{\ast}$,
so it is easy to calculate that
\begin{eqnarray*}
 v_{\ast}' = c\big[(am+cp+K(ar_{2}+cd))^{2}-4K(am+cp)(ar_{2}+cd-c r^{\ast}_{0})\big]^{-1/2},~
u_{\ast}' =m v_{\ast}'/c.
\end{eqnarray*}
Furthermore,
\begin{equation*}
   \tilde{A}_{1}=\left(\begin{array}{cc}
0  & \hat{a}_{12} \\
\hat{a}_{21} & \hat{a}_{22}  \\
\end{array}\right),
  \tilde{ B}_{1}=\left(\begin{array}{cc}
\hat{b}_{11}  & 0  \\
0 & 0   \\
\end{array}\right),
\end{equation*}
with $ \hat{a}_{12} =-\frac{K}{(1+K \tilde{v}_{\ast})^{3}}\big[(u_{\ast}'r_{0}^{\ast}+\tilde{u}_{\ast})(1+K\tilde{v}_{\ast})
 -2K\tilde{u}_{\ast}v_{\ast}'r_{0}^{\ast}\big]-pu_{\ast}'$,
   $\hat{a}_{21} = cv_{\ast}'$,
$\hat{a}_{22} =-m v_{\ast}'$,
$\hat{b}_{11}=-a u_{\ast}'$.
From (\ref{a916}), we get
\begin{eqnarray*}
F(0,0,U_{t})
=\left(\begin{array}{cc}
\tau^{\ast} \big[\tilde{\alpha}_{1}u_{t}(0)v_{t}(0)-au_{t}(0)u_{t}(-1)
+\tilde{\alpha}_{2}v_{t}^{2}(0)+\tilde{\alpha}_{3}v_{t}^{3}(0)+\tilde{\alpha}_{4}u_{t}(0)v_{t}^{2}(0) \big] \\
\tau^{\ast}\big[cu_{t}(0)v_{t}(0)-m v_{t}^{2}(0) \big]
\end{array}\right),\\
\end{eqnarray*}
where
\begin{eqnarray}
\begin{split}
\tilde{\alpha}_{1}=\frac{-Kr^{\ast}_{0}}{(1+K\tilde{v}_{\ast})^{2}}-p,
~\tilde{\alpha}_{2}=\frac{K^{2}r^{\ast}_{0}\tilde{u}_{\ast}}{(1+K\tilde{v}_{\ast})^{3}},~ \tilde{\alpha}_{3}=\frac{-K^{3}r^{\ast}_{0}\tilde{u}_{\ast}}{(1+K\tilde{v}_{\ast})^{4}},~
\tilde{\alpha}_{4}=\frac{K^{2}r^{\ast}_{0}}{(1+K\tilde{v}_{\ast})^{3}}.
\end{split}
\end{eqnarray}
Let $p_{1}(\theta) =(1,p_{12})^{T} e^{i\omega_{k_{1}}^{+}\tau^{\ast}\theta}$,
$p_{3}(\theta) =(1,p_{32})^{T} e^{i\omega_{k_{2}}^{-}\tau^{\ast}\theta}$ $(\theta\in[-1,0])$
be the eigenvectors of $\mathcal{A}_{m} (m=1,2)$ corresponding to the eigenvalue
$i\omega_{k_{1}}^{+}\tau^{\ast}$, $i\omega_{k_{2}}^{-}\tau^{\ast}$, respectively.
Besides, let $q_{1}(s) =D_{1}(1,q_{12}) e^{-i\omega_{k_{1}}^{+}\tau^{\ast}s}$,
$q_{3}(s) =D_{3}(1,q_{32})e^{-i\omega_{k_{2}}^{-}\tau^{\ast}s}$ $(s\in[0,1])$
be the eigenvectors of $\mathcal{A}_{m}^{\ast} (m=1,2)$ corresponding to the eigenvalue $i\omega_{k_{1}}^{+}\tau^{\ast}$, $i\omega_{k_{2}}^{-}\tau^{\ast}$, respectively. Via direct calculations, we can choose
$p_{12}=a_{21}^{\ast}/(i\omega_{k_{1}}^{+}+d_{2}\frac{k^{2}}{l^{2}}-a_{22}^{\ast})$,
$p_{32}=a_{21}^{\ast}/(i\omega_{k_{2}}^{-}+d_{2}\frac{k^{2}}{l^{2}}-a_{22}^{\ast})$,
$q_{12}=a_{12}^{\ast}/(i\omega_{k_{1}}^{+}+d_{2}\frac{k^{2}}{l^{2}}-a_{22}^{\ast})$,
$q_{32}=a_{12}^{\ast}/(i\omega_{k_{2}}^{-}+d_{2}\frac{k^{2}}{l^{2}}-a_{22}^{\ast})$,
$D_{1}=(1+p_{12}q_{12}+\tau^{\ast}b_{11}^{\ast}e^{-i\omega_{k_{1}}^{+}\tau^{\ast}})^{-1}$,
$D_{3}=(1+p_{32}q_{32}+\tau^{\ast}b_{11}^{\ast}e^{-i\omega_{k_{2}}^{-}\tau^{\ast}})^{-1}$.
By \citep{Du2018Double,Faria1995NormalBT}, we can easily show
that the usual normal form of the Hopf-Hopf bifurcation up to the
third order is as follows.
\begin{eqnarray}\label{f}
\begin{split}
&\dot{z}_{1} = i\omega_{k_{1}}^{+}\tau^{\ast}z_{1}+B_{11}\mu_{1}z_{1}+B_{21}\mu_{2}z_{1}+B_{2100}z_{1}^{2}z_{2}+B_{1011}z_{1}z_{3}z_{4},\\
&\dot{z}_{2} = -i\omega_{k_{1}}^{+}\tau^{\ast}z_{2}+\overline{B_{11}}\mu_{1}z_{2}+\overline{B_{21}}\mu_{2}z_{2}
+\overline{B_{2100}}z_{1}z_{2}^{2}+\overline{B_{1011}}z_{2}z_{3}z_{4},\\
&\dot{z}_{3} = i\omega_{k_{2}}^{-}\tau^{\ast}z_{3}+B_{13}\mu_{1}z_{3}+B_{23}\mu_{2}z_{3}+B_{0021}z_{3}^{2}z_{4}+B_{1110}z_{1}z_{2}z_{3},\\
& \dot{z}_{4} = -i\omega_{k_{2}}^{-}\tau^{\ast}z_{4}+\overline{B_{13}}\mu_{1}z_{4}+
\overline{B_{23}}\mu_{2}z_{4}+\overline{B_{0021}}z_{3}z_{4}^{2}+\overline{B_{1110}}z_{1}z_{2}z_{4},
\end{split}
\end{eqnarray}
where
$B_{11}=q_{1}(0)\Big[-\frac{k_{1}^{2}}{l^{2}}D p_{1}(0)+\tilde{A}_{0}p_{1}(0)+\tilde{B}_{0}p_{1}(-1)\Big]$,
$B_{21}=q_{1}(0)\tau^{\ast}\big[\tilde{A}_{1}p_{1}(0)+\tilde{B}_{1}p_{1}(-1)\big]$,
$B_{13}=q_{3}(0)\Big[-\frac{k_{2}^{2}}{l^{2}}D p_{3}(0)+\tilde{A}_{0}p_{3}(0)+\tilde{B}_{0}p_{3}(-1)\Big]$,
$B_{23}=q_{3}(0)\tau^{\ast}\big[\tilde{A}_{1}p_{3}(0)+\tilde{B}_{1}p_{3}(-1)\big]$.
The detailed derivation of the third-order normal form and expressions for $B_{2100}$, $B_{1011}$, $B_{0021}$, $B_{1110}$ are
given in appendix B.

To further analyze the bifurcation situation
in (\ref{f}), we use a cylindrical coordinate transformation
$z_{1}=\rho_{1}\cos\eta_{1}+i\rho_{1}\sin\eta_{1}$, $z_{2} = \rho_{1}\cos\eta_{1}-i\rho_{1}\sin\eta_{1}$,
$z_{3}=\rho_{2}\cos\eta_{2}+i\rho_{2}\sin\eta_{2}$, $z_{4} = \rho_{2}\cos\eta_{2}-i\rho_{2}\sin\eta_{2}$,
and define $\epsilon_{1}=\mathrm{Sign}(\mathrm{Re}B_{2100})$, $\epsilon_{2}=\mathrm{Sign}(\mathrm{Re}B_{0021})$,
and rescale $\bar{\rho}_{1}=\rho_{1}\sqrt{|B_{2100}|}$, $\bar{\rho}_{2}=\rho_{2}\sqrt{|B_{0021}|}$, $\bar{t}=t\epsilon_{1}$, then
drop the bars. System (\ref{f}) becomes,
\begin{eqnarray}\label{606a}
  \dot{\rho}_{1} =\rho_{1}(\upsilon_{1}+\rho_{1}^{2}+b_{0}\rho_{2}^{2}), ~
    \dot{\rho}_{2} = \rho_{2}(\upsilon_{2}+c_{0}\rho_{1}^{2}+d_{0}\rho_{2}^{2}),
\end{eqnarray}
where $\rho_{1},\rho_{2}>0$, and
\begin{eqnarray*}
\begin{split}
 & \upsilon_{1}=\epsilon_{1}(\mathrm{Re} B_{11}\mu_{1}+\mathrm{Re} B_{21}\mu_{2})=\epsilon_{1}\Big[\mathrm{Re} B_{11}(\tau-\tau^{\ast})+\mathrm{Re} B_{21}(r_{0}-r_{0}^{\ast})\Big], \\
 & \upsilon_{2}=\epsilon_{1}(\mathrm{Re} B_{13}\mu_{1}+\mathrm{Re} B_{23}\mu_{2})=\epsilon_{1}\Big[\mathrm{Re} B_{13}(\tau-\tau^{\ast})+\mathrm{Re} B_{23}(r_{0}-r_{0}^{\ast})\Big], \\
 &b_{0}=\frac{\epsilon_{1}\epsilon_{2}\mathrm{Re}B_{1011}}{\mathrm{Re}B_{0021}}, c_{0}=\frac{\mathrm{Re}B_{1110}}{\mathrm{Re}B_{2100}},
d_{0}=\epsilon_{1}\epsilon_{2}=\pm1.
\end{split}
\end{eqnarray*}
There are twelve distinct types of unfolding
for (\ref{606a}) due to the different signs of $d_{0}$, $b_{0}$, $c_{0}$ and $d_{0}-b_{0}c_{0}$.
The readers may
find the detailed phase portraits in chapter 7 in the original book \citep{J1983Nonlinear}.

\section{Numerical examples}\label{a1112}
In this section, we show some numerical simulations to support the theoretical results obtained in sections \ref{105d} and \ref{a1019}. The rich dynamics such as periodic solutions near Hopf or Hopf-Hopf bifurcations are illustrated. Moreover, we show the coexistence of stationary solution and periodic solution near a Bautin bifurcation point.
\subsection{Simulations about Hopf bifurcation}
\begin{figure}[!htbp]
\centering{\includegraphics[width=0.8\textwidth]{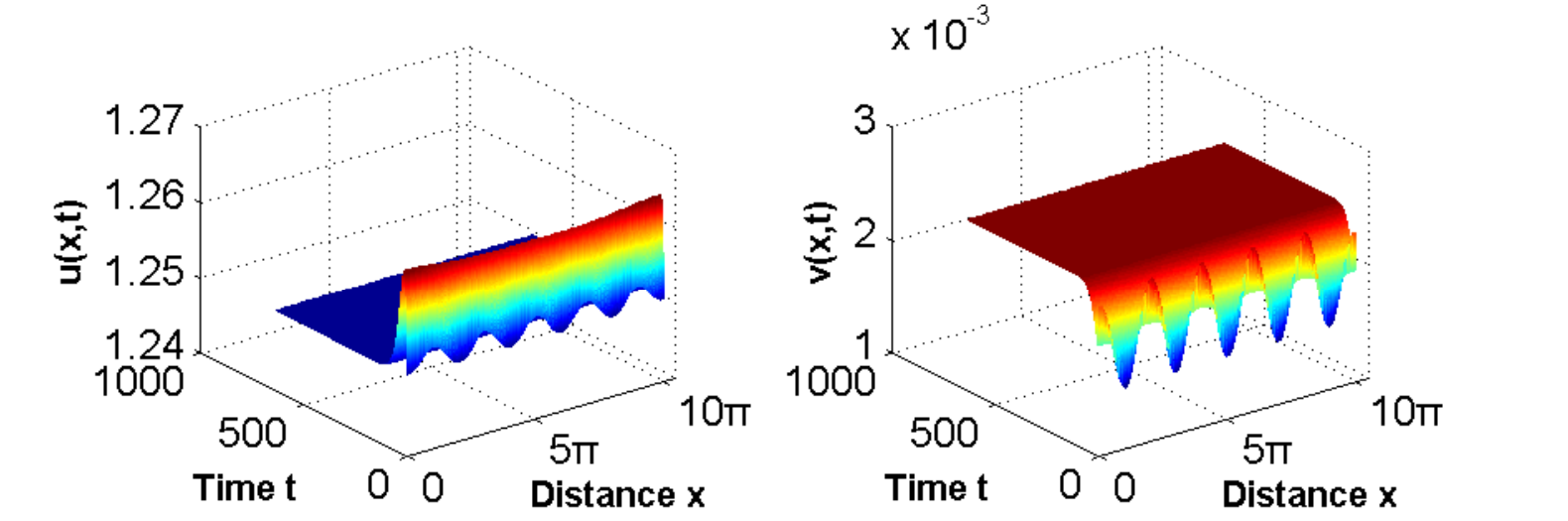}}
\caption{The positive equilibrium is locally asymptotically stable, where the initial functions are $(u_{0}(x,t),v_{0}(x,t))=(1.25+0.001\mathrm{cos}x,0.002+0.001\mathrm{cos}x)$ and $\tau=4<\bar{\tau}$.}
\label{b929}
\end{figure}

\begin{figure}[!htbp]
\centering{\includegraphics[width=0.8\textwidth]{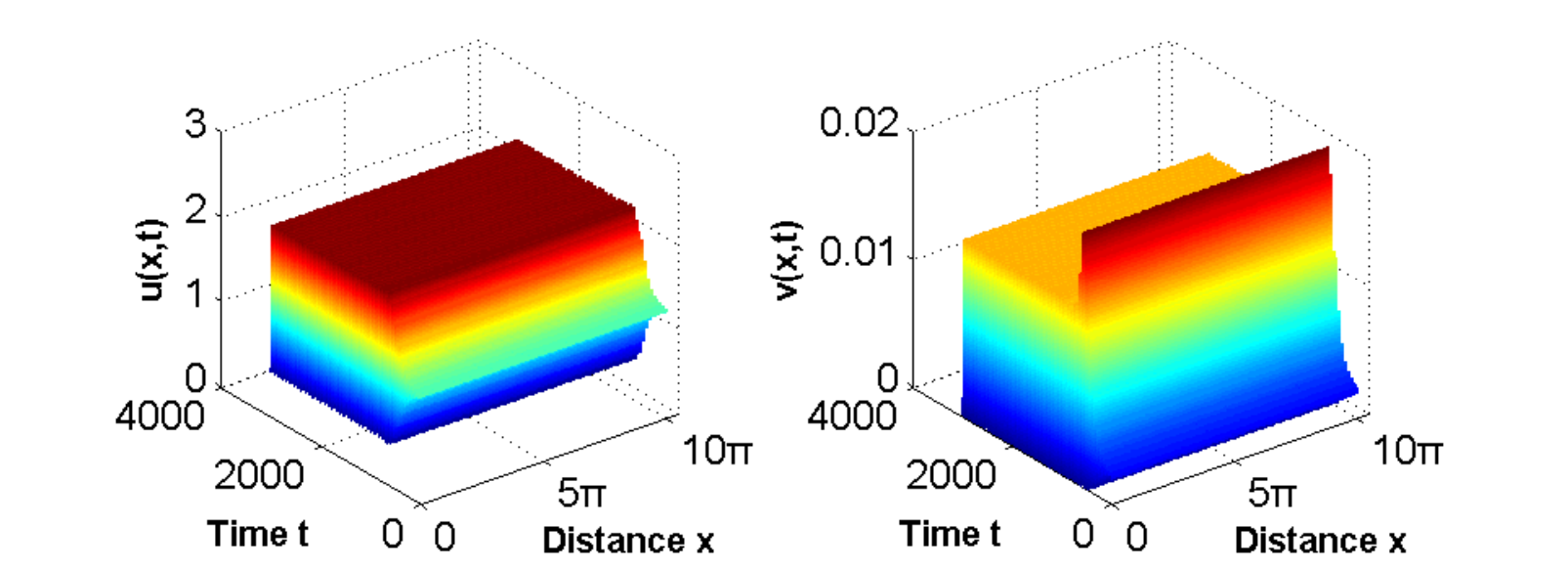}}
\caption{The bifurcating periodic solutions are stable, where the initial functions are $(u_{0}(x,t),v_{0}(x,t))=(1.25+0.001\mathrm{cos}x,0.002+0.001\mathrm{cos}x)$ and $\tau=19>\bar{\tau}$.}
\label{c929}
\end{figure}

Now we carry out some simulations for system (\ref{a}). Fix parameter $r_{0}=0.12$ and choose
\begin{eqnarray}\label{929a}
d_{1}=0.3,d_{2}=0.5,K=10,d=0.04, a=0.06, p=0.8, r_{2}=0.5, c=0.4, m=0.1, l=10.
\end{eqnarray}
For this set of parameter values we find that (H3) holds and the positive equilibrium is $E_{\ast}(1.2506, 0.0025)$. By direct calculation, we find
$\bar{\tau}=\tau_{0}^{0+}\approx15.7797$, $\omega_{0}^{+}=0.0996$, $c_{1}(0)\approx -0.0022-0.0032i$, $\mu_{2}\approx0.00001$ and $\beta_{2}=-0.0045$.
By Theorem \ref{a912}, if $\tau\in[0,\bar{\tau})$, then the equilibrium $E_{\ast}$ is locally asymptotically stable (see Fig. \ref{b929}).
By Theorem \ref{a1113}, the direction of the Hopf
bifurcation is forward when $\tau=\bar{\tau}$, and the bifurcating periodic solutions are orbitally  asymptotically stable when $\tau>\bar{\tau}$ (see Fig. \ref{c929}).
Besides, if we choose $r_{0}=1$, we can find $E_{\ast}$ is asymptotically stable when
$\tau\in[0,3.6892)\cup(10.4124,16.4128)\cup(25.4181,29.1364)\cup(40.4235,41.8598)$ by
applying the results of Theorem \ref{a912} (iii). These stable regions are shown in Fig. \ref{708a}.
\subsection{Simulations about Hopf-Hopf bifurcation}
\begin{figure}[!htbp]
\centering{\includegraphics[width=0.7\textwidth]{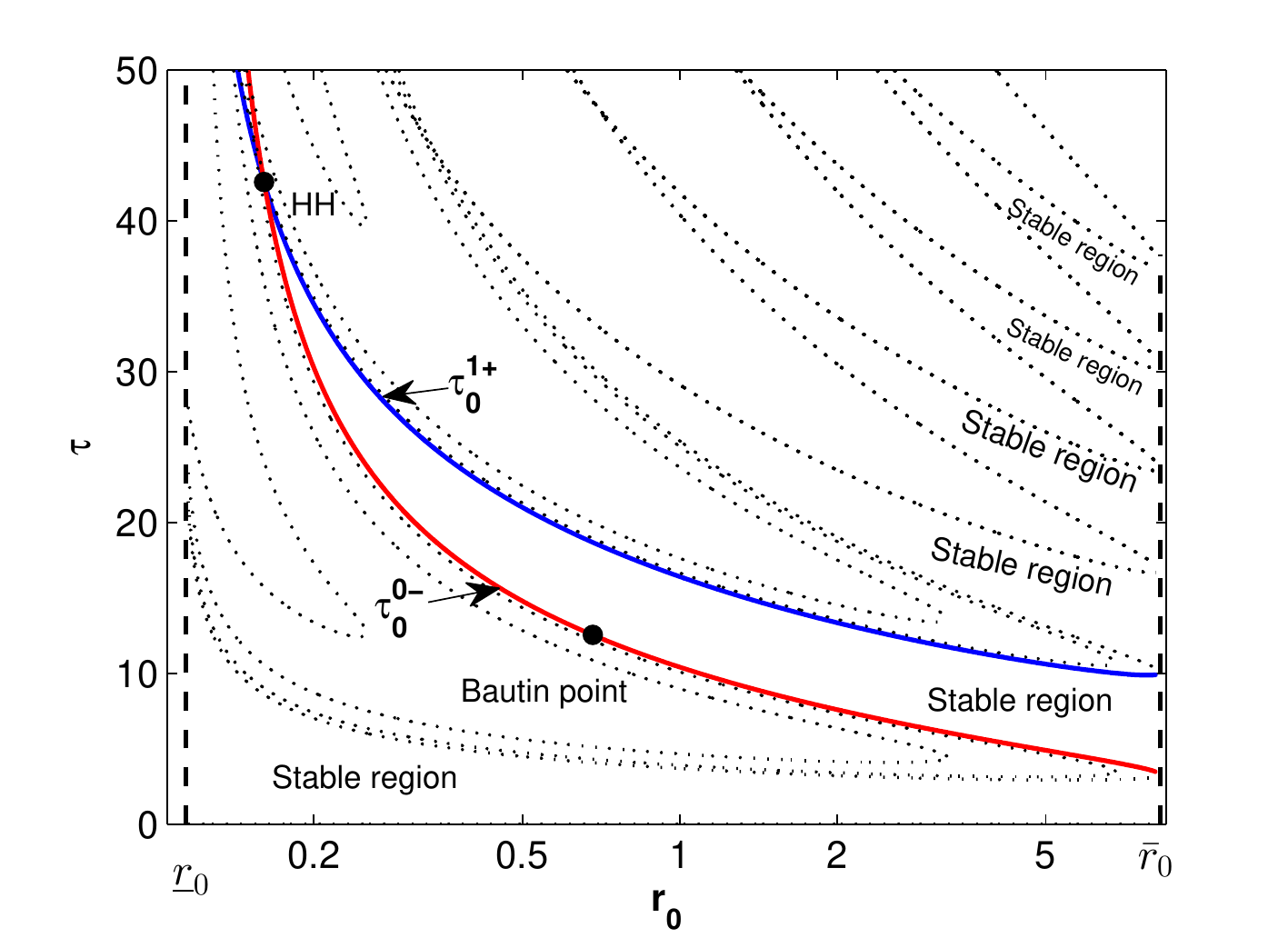}}
\caption{The partial bifurcation set on the $\tau-r_{0}$ plane. All Hopf bifurcation curves are marked by dotted curves, except the blue curve and red curve stand for $\tau_{0}^{1+}$ and $\tau_{0}^{0-}$, which intersect at the double Hopf-Hopf bifurcation point HH.
}
\label{708a}
\end{figure}

In this subsection, we want to show the rich dynamics such as periodic and quasi-periodic oscillations near Hopf-Hopf bifurcation
in system ($\ref{a}$) by a group of simulations. We regard $\tau$ and $r_{0}$ as bifurcation parameters, other parameters are the same as in (\ref{929a}),
the bifurcation diagram is shown in Fig. \ref{708a}. Two boundary points can be easily calculated as $\underline{r}_{0}=0.115$, $\bar{r}_{0}=8.2716$ for $k=0$, i.e. $\underline{r}_{0}<r_{0}<\bar{r}_{0}$, and Hopf bifurcation may occur between them.
We can also find $\mathrm{(H3)}$ holds for all $0\leq k\leq3$ when $0.115<r_{0}<0.2519$.
When $r_{0}^{\ast}=0.1606$, $\tau_{0}^{0-}$ firstly intersects with $\tau_{0}^{1+}$ at the Hopf-Hopf bifurcation point HH when $k=0$. For HH,
we obtain $\omega_{0}^{+}=0.1848$, $\omega_{0}^{-}=0.1095$ and $\tau^{\ast}=\tau_{0}^{0-}=\tau_{0}^{1+}=42.5794$.

\begin{figure}[!htbp]
\centering{\includegraphics[width=0.7\textwidth]{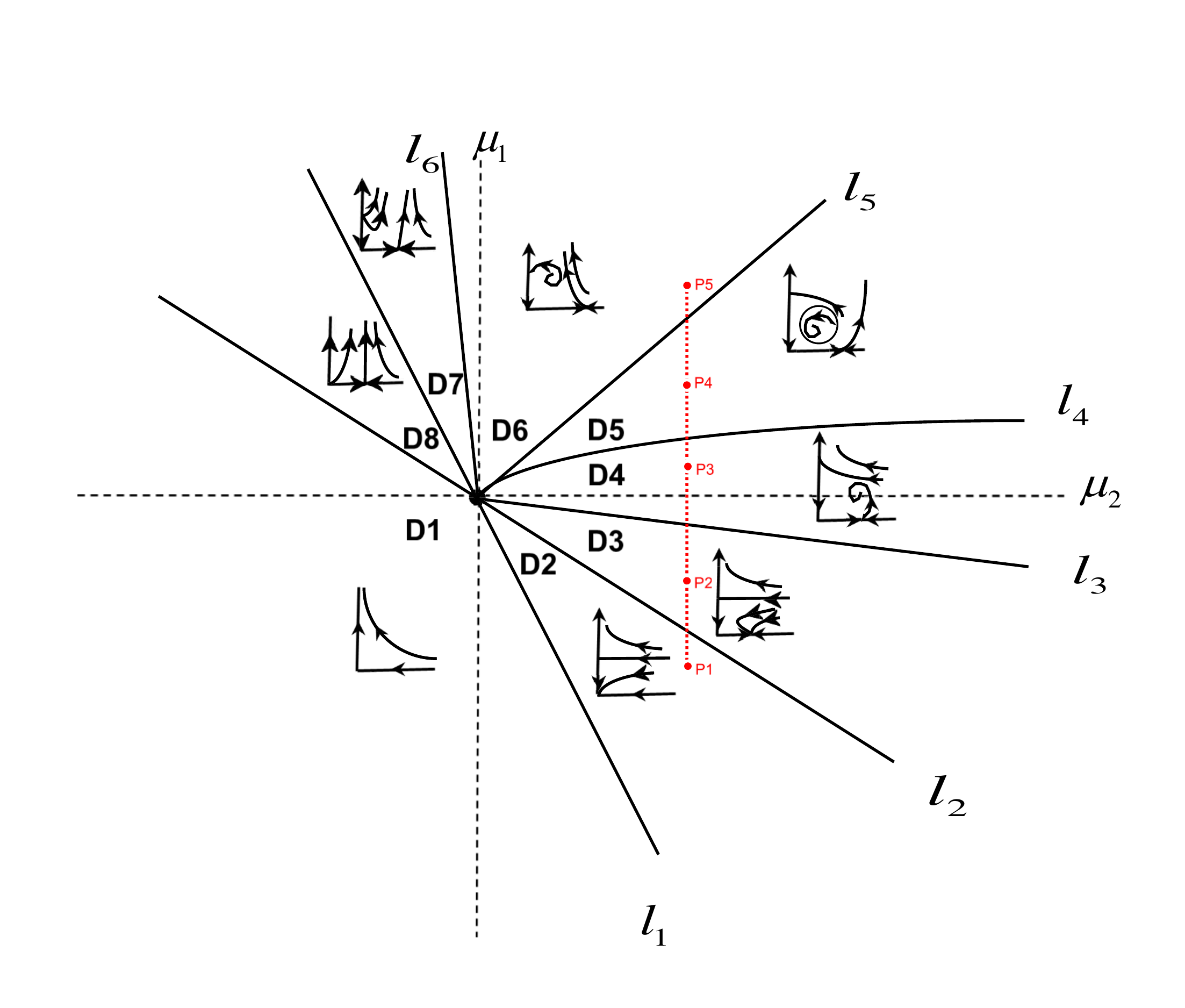}}
\caption{Complete bifurcation sets near the Hopf-Hopf point HH for system (\ref{a}).}
\label{b1021}
\end{figure}

\begin{figure}[!htbp]
\centering{\includegraphics[width=0.8\textwidth]{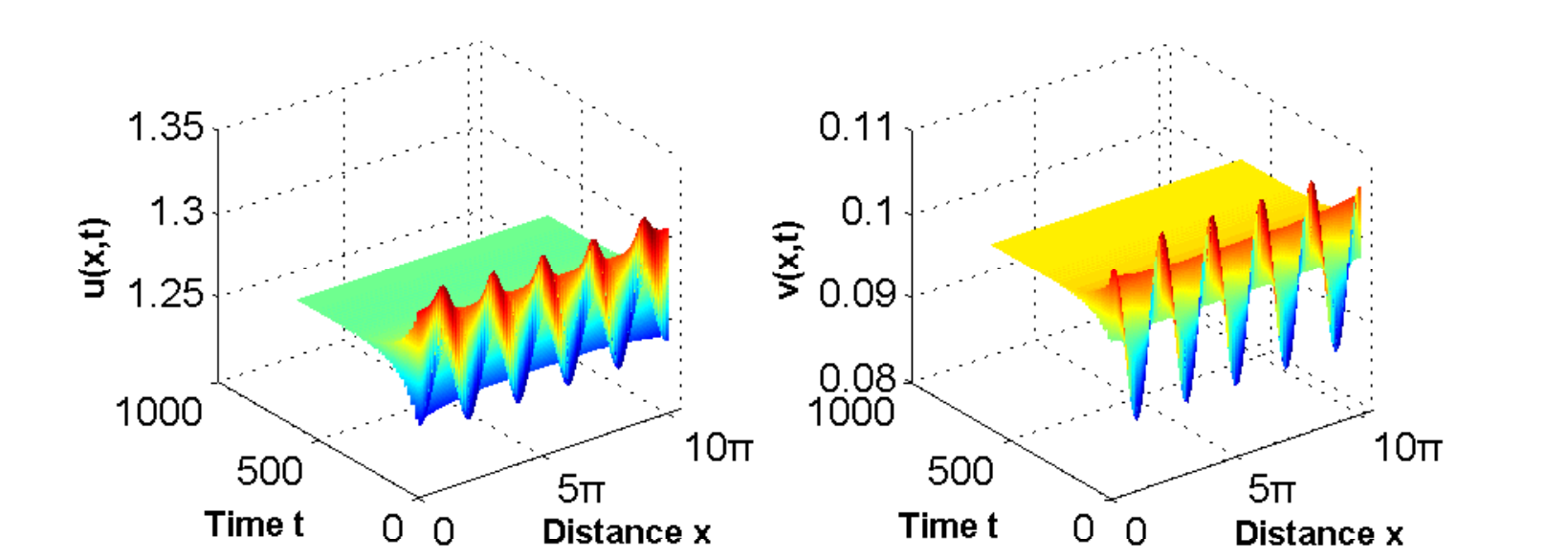}}
\caption{When $\tau=4$, $r_{0}=0.4$ (see P1 in Fig. \ref{b1021}),
the positive constant stationary solution $E_{\ast}$ of system (\ref{a}) is locally asymptotically stable,
the initial functions are $(u_{0}(x,t),v_{0}(x,t))=(1.25+0.02\mathrm{cos}x,0.1+0.02\mathrm{cos}x)$.}
\label{629d}
\end{figure}

\begin{figure}[!htbp]
\centering{\includegraphics[width=0.8\textwidth]{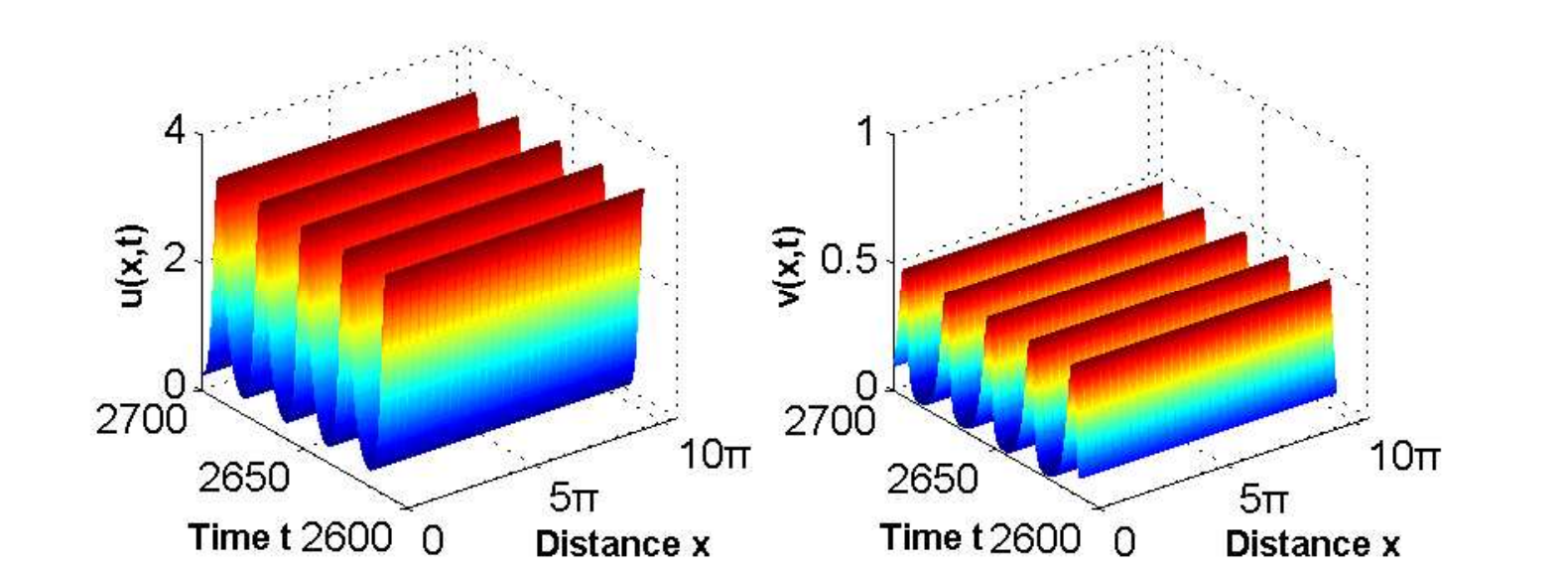}}
\caption{When $\tau=25$, $r_{0}=0.4$ (see P2 in Fig. \ref{b1021}), spatially homogeneous periodic solutions of system (\ref{a}) are stable,
the initial functions are $(u_{0}(x,t),v_{0}(x,t))=(1.25+0.02\mathrm{cos}x,0.1+0.02\mathrm{cos}x)$.}
\label{aa924}
\end{figure}

\begin{figure}[!htbp]
\centering{\includegraphics[width=0.8\textwidth]{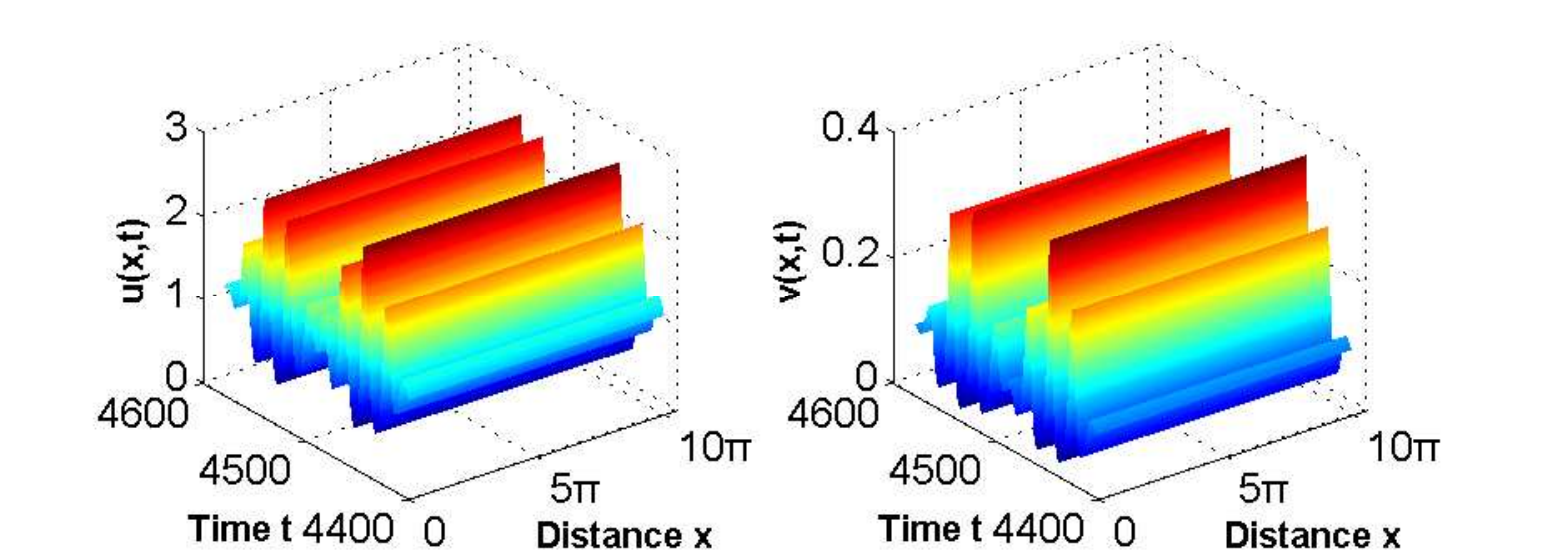}}
\caption{When $\tau=42.4$, $r_{0}=0.4$ (see P3 in Fig. \ref{b1021}), spatially homogeneous quasi-periodic solutions of system (\ref{a}) are unstable,
the initial functions are $(u_{0}(x,t),v_{0}(x,t))=(1.25+0.02\mathrm{cos}x,0.1+0.02\mathrm{cos}x)$.}
\label{bb924}
\end{figure}

\begin{figure}[!htbp]
\centering{\includegraphics[width=1\textwidth]{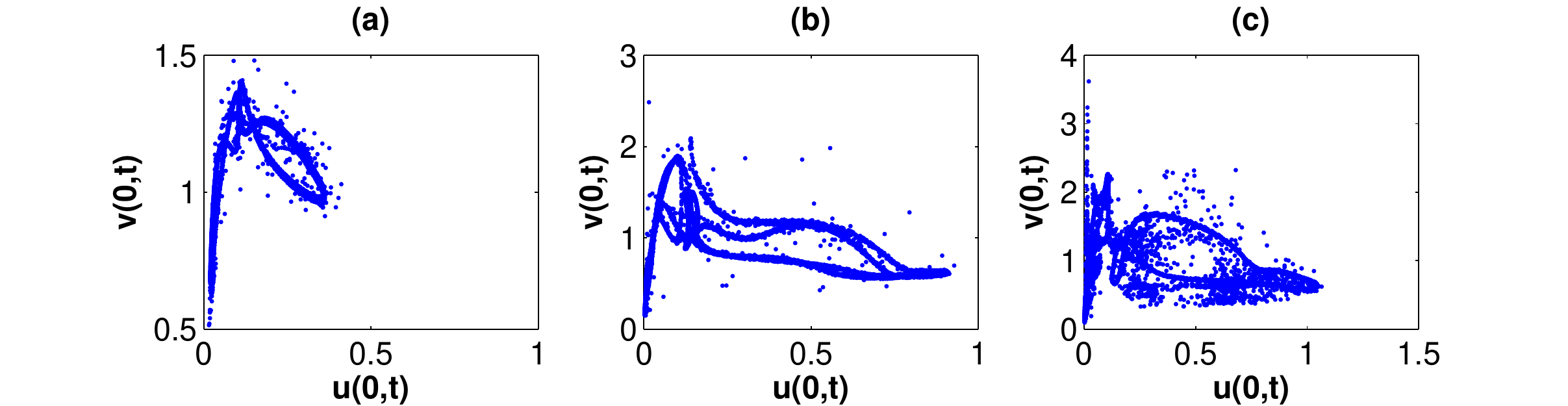}}
\caption{(a) There is a quasi-periodic solution of system (\ref{a}) when $\tau=42.4$, $r_{0}=0.4$ (see P3);
(b) There is a quasi-periodic solution of system (\ref{a}) when $\tau=46.6$, $r_{0}=0.4$ (see P4);
(c) There is a strange attractor of system (\ref{a}) when $\tau=47.4$, $r_{0}=0.4$ (see P5).
The initial values are all $(u_{0}(x,t),v_{0}(x,t))=(1.25+0.02\mathrm{cos}x,0.1+0.02\mathrm{cos}x)$.}
\label{cc924}
\end{figure}

After some straightforward calculations for system (\ref{f}), we obtain
$B_{11}=0.07218 + 0.03765i$, $B_{21}=22.59131+12.20126i$,
$B_{13}=-0.05265 + 0.0454i$, $B_{23}=-28.24841+23.62431i$,
$B_{2100}=-0.07041-0.04591i$, $B_{1011}=-0.34767- 0.32768i$,
$B_{0021}=0.1865-0.04057i$, $B_{1110}=0.10178- 0.24012i$.
and
\begin{eqnarray}\label{y614}
\epsilon_{1}=-1,~\epsilon_{2}=1,~d_{0}=-1,~b_{0}=1.8642, c_{0}=-1.4456, d_{0}-b_{0}c_{0}=1.6949.
\end{eqnarray}
It follows from (\ref{y614}) that the case VIa occurs in chapter 7 in \citep{J1983Nonlinear}.
We consider the dynamics of system (\ref{a}) when $(\mu_{1},\mu_{2})=(\tau-\tau^{\ast}, r_{0}-r_{0}^{\ast})$ is sufficiently  close to $(0,0)$ and divide the region $\mathbf{D}=\{(\mu_{1},\mu_{2})\}$ into eight parts by lines $l_{1}$, $l_{2}$, curve $l_{4}$, and half-lines $l_{3}$, $l_{5}$, $l_{6}$:
\begin{eqnarray*}
&& l_{1}: \mu_{1}=-536.532\mu_{2},~~l_{2}:\mu_{1}=-312.9857\mu_{2},\\
&& l_{3}: \mu_{1}=-85.2815\mu_{2} ~(\mu_{2}\geq0),
~l_{4}: \mu_{1}=997.7215\mu_{2}+o(\mu_{2})~(\mu_{2}\geq0),\\
&& l_{5}: \mu_{1}=997.7215\mu_{2} ~(\mu_{2}\geq0),~l_{6}:\mu_{1}=-1157.8669\mu_{2}~(\mu_{2}\geq0).
\end{eqnarray*}
These regions are illustrated in Fig. \ref{b1021}. Then we can obtain the following information for system (\ref{a}).
\begin{enumerate}
\item[$(\mathbf{a})$] In region $\mathrm{D2}$, the trivial equilibrium is a sink and corresponds to the positive constant stationary solution $E_{\ast}$ of system (\ref{a}) (see Fig. \ref{629d}).
\item[$(\mathbf{b})$] A stable periodic solution appears as $(\mu_{1},\mu_{2})$ crosses $l_{2}$ from $\mathrm{D2}$ to $\mathrm{D3}$ via supercritical Hopf bifurcations.
    We simulate the spatially homogeneous stable periodic solutions when $P_{2}$ is chosen in D3 (see Fig. \ref{aa924}).
\item[$(\mathbf{c})$] There is a quasi-periodic solution on the two-dimensional torus when $P_{3}$ is chosen in $\mathrm{D4}$. (see Fig. \ref{bb924}).
    Since the Hopf-Hopf bifurcation occurs at the intersection of two Hopf bifurcation curves both with wave number $k=0$, the periodic solutions of system (\ref{a}) near the bifurcation point are always spatially homogeneous.
Thus, the solution curves of $(u(0,t), v(0,t))$ can represent the dynamical behavior of the whole solution $(u(x,t), v(x,t))$ (see \citep{Du2018Double}).
On the Poincar$\acute{\mathrm{e}}$ section $u(0,t-\tau)=u_{\ast}$, where $u_{\ast}$ is the constant steady solution in (\ref{b}). If we further perturb the parameters,
we find that the points on the Poincar$\acute{\mathrm{e}}$ section exhibit
quasi-periodic behavior (see Fig. \ref{cc924}a,b). Besides,
there exists a strange attractor of system (\ref{a}), see Fig. \ref{cc924}c.
\end{enumerate}
\subsection{Simulations about Bautin bifurcation}
\begin{figure}[!htbp]
\centering{a\includegraphics[width=0.8\textwidth]{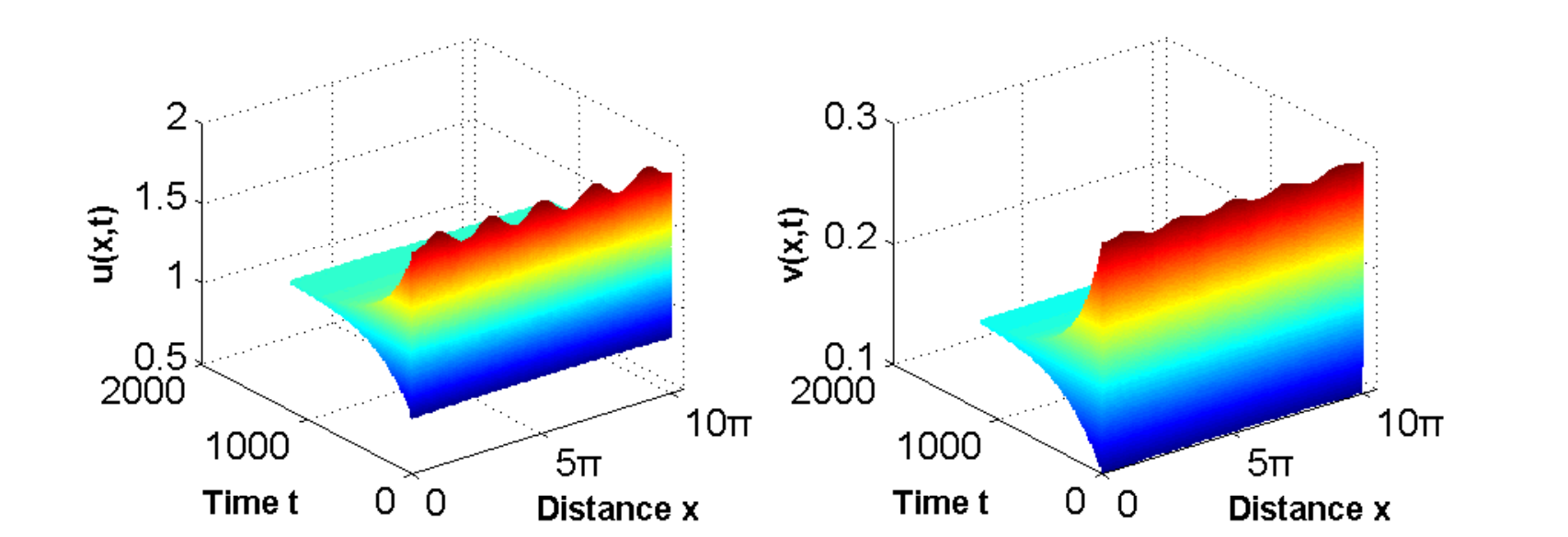}}
{b\includegraphics[width=0.8\textwidth]{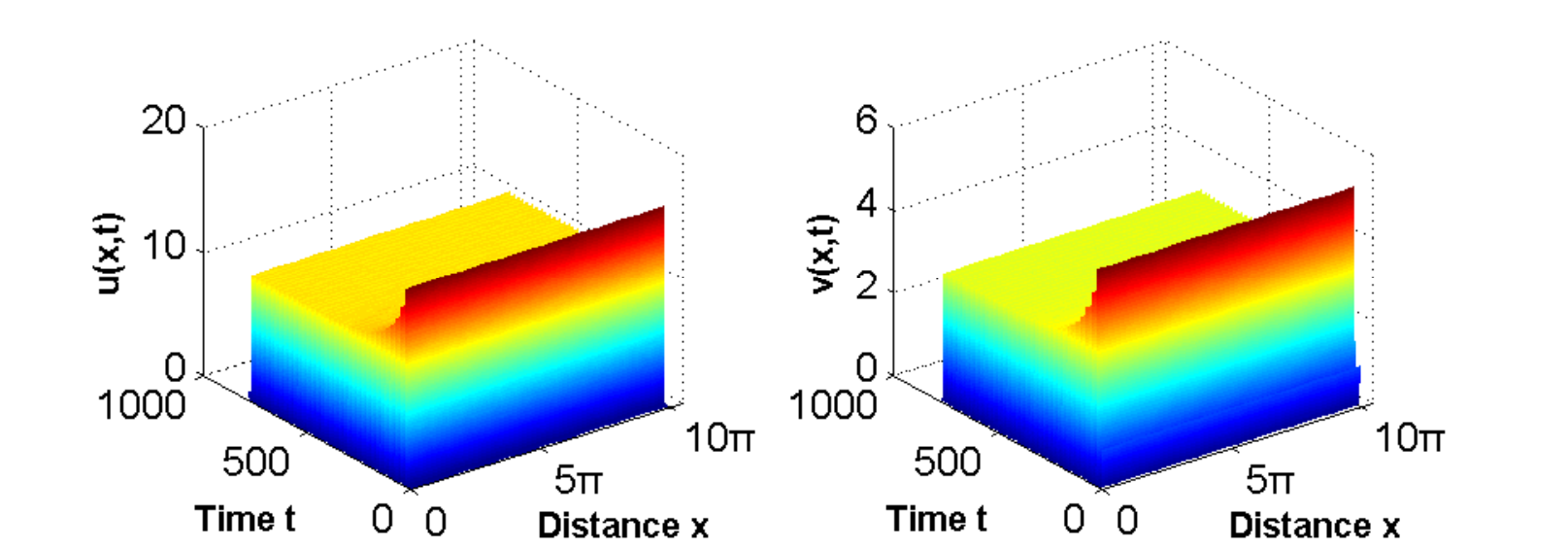}}
\caption{The stable constant stationary solution $E_{\ast}$ and spatially homogeneous stable periodic solution of system (\ref{a}) with different initial values coexist.
(a) The initial functions are $(u_{0}(x,t),v_{0}(x,t))=(1.25+0.02\mathrm{cos}x,0.1+0.02\mathrm{cos}x)$; (b) The initial functions are $(u_{0}(x,t),v_{0}(x,t))=(6+0.02\mathrm{cos}x,3+0.02\mathrm{cos}x)$.
Here  $\tau=13$, $r_{0}=0.7$, and all the other parameters are given by (\ref{929a}).}
\label{a1021}
\end{figure}

The red Hopf bifurcation curve below HH is subcritical near the Hopf-Hopf bifurcation point, see Fig. \ref{708a}. We find when $(r_{0},\tau)=(0.682,12.545)$, $\mathrm{Re}(c_{1}(0))=0$ holds, i.e.,
system (\ref{a}) may undergo a Bautin bifurcation at $(0.682,12.545)$.
In fact, we find that near the Bautin bifurcation point,
the stable constant stationary solution $E_{\ast}$ and spatially homogeneous stable periodic solution of system (\ref{a}) with different initial values coexist, see Fig. \ref{a1021}.
\section{Conclusions and Discussions}\label{b1112}
In this paper, we mainly study the Hopf bifurcation and nonresonant Hopf-Hopf bifurcation in a predator-prey system with fear effect in prey.
We compute the normal form near the Hopf and Hopf-Hopf bifurcation points and give an explicit algorithm for calculating the four key variables at the Hopf-Hopf bifurcation: $b_{0}, c_{0}, d_{0}$ and $d_{0}-b_{0}c_{0}$.
Detailed dynamics near the critical point are obtained by drawing the corresponding bifurcation set.
Through bifurcation analysis, we find such a predator-prey system with fear effect has very rich dynamics, including periodic and quasi-periodic oscillations.
With the aid of bifurcation set, we also find that a strange attractor appears while the
quasi-periodic oscillation on three-torus vanishes through a
saddle connection bifurcation.

Near a Bautin bifurcation point, we numerically find the coexistence of stable steady state and stable periodic oscillation, which means solutions with different initial functions converge to a stable
constant stationary solution or spatially homogeneous periodic solution, respectively.
In fact, the detailed algorithm for deriving the normal form of Bautin bifurcation in the reaction-diffusion system with time delay could be given,
which will be left as a future work.

\section*{Acknowledgement}
This research is supported by National Natural Science Foundation of China (No.11701120 and No.11771109).

\appendix
\section{Properties of Hopf bifurcation}
We can rewrite system (\ref{a916}) as follows:
\begin{eqnarray}\label{a917}
\dot{U}(t)=(\bar{\tau}+\mu)D\Delta U(t)+(\bar{\tau}+\mu)(\tilde{A}U(t)+\tilde{B}U(t-1))+F(\mu,U_{t}),
\end{eqnarray}
where
\begin{eqnarray*}
F(\mu,\phi)=
(\bar{\tau}+\mu)\left(\begin{array}{cc}
\alpha_{1}\phi_{1}(0)\phi_{2}(0)-a\phi_{1}(0)\phi_{1}(-1)+\alpha_{2}\phi_{2}^{2}(0)
+\alpha_{3}\phi_{2}^{3}(0)+\alpha_{4}\phi_{1}(0)\phi_{2}^{2}(0)\\
c\phi_{1}(0)\phi_{2}(0)-m\phi_{2}^{2}(0)
\end{array}\right)
\end{eqnarray*}
for $(\phi_{1},\phi_{2})\in\mathcal{C}$. Then we decompose the space $\mathcal{C}$ as $C=P\oplus Q$, where
$P=\{zp\gamma_{n_{0}}(x)+\bar{z}\bar{p}\gamma_{n_{0}}(x)|z\in\mathbb{C}\}$,
$Q=\{\phi\in\mathcal{C}|(q\gamma_{n_{0}}(x),\phi)=0 ~\mathrm{and}~ (\bar{q}\gamma_{n_{0}}(x),\phi)=0\}$.
Thus, system (\ref{a917}) could be rewritten as
$U_{t}=z(t)p(\cdot)\gamma_{n_{0}}(x)+\bar{z}(t)\bar{p}(\cdot)\gamma_{n_{0}}(x)+w(t,\cdot)$,
where
\begin{eqnarray}\label{a922}
z(t)=(q\gamma_{n_{0}}(x),U_{t}),~w(t,\theta)=U_{t}(\theta)-2\mathrm{Re}\{z(t)p(\theta)\gamma_{n_{0}}(x)\},
\end{eqnarray}
then we have $\dot{z}(t)=i\omega_{n_{0}}\bar{\tau}z(t)+\bar{q}(0)\langle F(0,U_{t}),\beta_{n_{0}}\rangle$.
There exists a center manifold $\mathcal{C}_{0}$ and we can write $w$ near $(0,0)$ as follows.
\begin{eqnarray}\label{c922}
w(t,\theta)=w(z(t),\bar{z}(t),\theta)=w_{20}(\theta)\frac{z^{2}}{2}+w_{11}(\theta)z\bar{z}+w_{02}(\theta)\frac{\bar{z}^{2}}{2}+\cdots,
\end{eqnarray}
the system restricted to the center manifold is given by
$\dot{z}(t)=i\omega_{n_{0}}\bar{\tau}z(t)+g(z,\bar{z})$,
and denote $g(z,\bar{z})=g_{20}\frac{z^{2}}{2}+g_{11}z\bar{z}
+g_{02}\frac{\bar{z}^{2}}{2}+g_{21}\frac{z^{2}\bar{z}}{2}+\cdots$.
By direct calculation, we obtain
\begin{eqnarray*}
&&g_{20}=2\bar{\tau}M(\alpha_{1}p_{1}-ae^{-i\omega_{n_{0}}\bar{\tau}}+\alpha_{2}p_{1}^{2}+c q_{2}p_{1}-m q_{2}p_{1}^{2})\int^{l\pi}_{0}\gamma^{3}_{n_{0}}(x)dx,\\
&&g_{11}=\bar{\tau}M\big[(\alpha_{1}+c q_{2})(p_{1}+\bar{p}_{1})-a(e^{i\omega_{n_{0}}\bar{\tau}}
+e^{-i\omega_{n_{0}}\bar{\tau}})+2p_{1}\bar{p}_{1}(\alpha_{2}-m q_{2})\big]\int^{l\pi}_{0}\gamma^{3}_{n_{0}}(x)dx,\\
&&g_{02}=2\bar{\tau}M(\alpha_{1}\bar{p}_{1}-ae^{i\omega_{n_{0}}\bar{\tau}}+\alpha_{2}\bar{p}_{1}^{2}+c q_{2}\bar{p}_{1}-m q_{2}\bar{p}_{1}^{2})\int^{l\pi}_{0}\gamma^{3}_{n_{0}}(x)dx,\\
&&g_{21}=2\bar{\tau}M\bigg[(3\alpha_{3}p_{1}^{2}\bar{p}_{1}+\alpha_{4}(2p_{1}\bar{p}_{1}
+p_{1}^{2}))\int^{l\pi}_{0}\gamma^{4}_{n_{0}}(x)dx+Q_{1}\int^{l\pi}_{0}\gamma^{2}_{n_{0}}(x)dx\bigg],
\end{eqnarray*}
where
\begin{eqnarray*}
Q_{1}&=&(\alpha_{1}+c q_{2})\bigg[w_{11}^{(2)}(0)+\frac{w_{20}^{(2)}(0)}{2}
+\frac{w_{20}^{(1)}(0)\bar{p}_{1}}{2}+p_{1}w_{11}^{(1)}(0)\bigg]
+(\alpha_{2}-m q_{2})\big[2p_{1}w_{11}^{(2)}(0)+\bar{p}_{1}w_{20}^{(2)}(0)\big]\\
&-&a\bigg[w_{11}^{(1)}(-1)+\frac{w_{20}^{(1)}(-1)}{2}+\frac{w_{20}^{(1)}(0)}{2}
e^{i\omega_{n_{0}}\bar{\tau}}+w_{11}^{(1)}(0)e^{-i\omega_{n_{0}}\bar{\tau}}\bigg],
\end{eqnarray*}
In order to get $g_{21}$, we need to compute $w_{20}$ and $w_{11}$. From (\ref{a922}), we
have
\begin{eqnarray}\label{b922}
\dot{w}=\dot{U}_{t}-\dot{z}p\gamma_{n_{0}}(x)-\dot{\bar{z}}\bar{p}\gamma_{n_{0}}(x)
\dot{=}\mathcal{A}w+H(z,\bar{z},\theta),
\end{eqnarray}
where
$H(z,\bar{z},\theta)=H_{20}(\theta)\frac{z^{2}}{2}
+H_{11}(\theta)z\bar{z}+H_{02}(\theta)\frac{\bar{z}^{2}}{2}+\cdots$.
Comparing the coefficients of (\ref{b922}) with (\ref{c922}), we obtain
\begin{eqnarray}\label{d922}
(\mathcal{A}-2i\omega_{n_{0}}\bar{\tau} I)w_{20}(\theta)=-H_{20}(\theta),~\mathcal{A} w_{11}(\theta)=-H_{11}(\theta),\cdots.
\end{eqnarray}
By (\ref{d922}), we have
\begin{eqnarray*}
&&w_{20}(\theta)=\frac{-g_{20}}{i\omega_{n_{0}}\bar{\tau}}p(0)
e^{i\omega_{n_{0}}\bar{\tau}\theta}\gamma_{n_{0}}(x)-
\frac{\bar{g}_{02}}{3i\omega_{n_{0}}\bar{\tau}}\bar{p}(0)
e^{-i\omega_{n_{0}}\bar{\tau}\theta}\gamma_{n_{0}}(x)+
E_{1}e^{2i\omega_{n_{0}}\bar{\tau}\theta},\\
&&w_{11}(\theta)=\frac{g_{11}}{i\omega_{n_{0}}\bar{\tau}}p(0)
e^{i\omega_{n_{0}}\bar{\tau}\theta}\gamma_{n_{0}}(x)-
\frac{\bar{g}_{11}}{i\omega_{n_{0}}\bar{\tau}}\bar{p}(0)
e^{-i\omega_{n_{0}}\bar{\tau}\theta}\gamma_{n_{0}}(x)+E_{2},
\end{eqnarray*}
Denote
$E_{1}=\sum^{\infty}_{n=0}E_{1}^{n}\gamma_{n}(x)$, $E_{2}=\sum^{\infty}_{n=0}E_{2}^{n}\gamma_{n}(x)$,
$E_{1}^{n}$ and $E_{2}^{n}$ could be calculated by
\begin{eqnarray*}
&&E_{1}^{n}=\bigg(2i\omega_{n_{0}}\bar{\tau}I-\int_{-1}^{0}
e^{2i\omega_{n_{0}}\bar{\tau}\theta}\mathrm{d}\eta^{n_{0}}(\theta,\bar{\tau})\bigg)^{-1}\langle \tilde{F}_{20},\beta_{n}\rangle,\\
&&E_{2}^{n}=-\bigg(\int_{-1}^{0}\mathrm{d}\eta^{n_{0}}(\theta,\bar{\tau})\bigg)^{-1}\langle \tilde{F}_{11},\beta_{n}\rangle,~n=0,1,\cdots\\
\end{eqnarray*}
where
\begin{eqnarray*}
\langle \tilde{F}_{20},\beta_{n}\rangle=
\begin{cases}
\frac{1}{\sqrt{l\pi}}\hat{F}_{20},&n_{0}\neq0,n=0,\\
\frac{1}{\sqrt{2l\pi}}\hat{F}_{20},&n_{0}\neq0,n=2n_{0},\\
\frac{1}{\sqrt{l\pi}}\hat{F}_{20},&n_{0}=0,n=0,\\
0,&\mathrm{other},
\end{cases}
\langle \tilde{F}_{11},\beta_{n}\rangle=
\begin{cases}
\frac{1}{\sqrt{l\pi}}\hat{F}_{11},&n_{0}\neq0,n=0,\\
\frac{1}{\sqrt{2l\pi}}\hat{F}_{11},&n_{0}\neq0,n=2n_{0},\\
\frac{1}{\sqrt{l\pi}}\hat{F}_{11},&n_{0}=0,n=0,\\
0,&\mathrm{other},
\end{cases}
\end{eqnarray*}
and
\begin{eqnarray*}
\hat{F}_{20}=2\left(\begin{array}{cc}
\alpha_{1}p_{1}-ae^{-i\omega_{n_{0}}\bar{\tau}}+\alpha_{2}p_{1}^{2}\\
c p_{1}-m p_{1}^{2}\\
\end{array}\right),~
\hat{F}_{11}=\left(\begin{array}{cc}
\alpha_{1}(p_{1}+\bar{p}_{1})-a(e^{i\omega_{n_{0}}\bar{\tau}}
+e^{-i\omega_{n_{0}}\bar{\tau}})+2\alpha_{2}p_{1}\bar{p}_{1}\\
c (p_{1}+\bar{p}_{1})-2m p_{1}\bar{p}_{1}\\
\end{array}\right).
\end{eqnarray*}
Therefore, $g_{21}$ could be represented explicitly.

\section{Derivation of third-order normal form of Hopf-Hopf bifurcation}
For system (\ref{f}) we have (see \citep{Du2018Double}),
\begin{eqnarray*}
&&B_{2100}= C_{2100}+\frac{3}{2}( D_{2100}+E_{2100}),~B_{1011}= C_{1011}+\frac{3}{2}(D_{1011}+E_{1011}),\\
&&B_{0021}=C_{0021}+\frac{3}{2}( D_{0021}+E_{0021}),~B_{1110}= C_{1110}+\frac{3}{2}(D_{1110}+E_{1110}).
\end{eqnarray*}
Firstly, by some direct calculation we obtain
\begin{eqnarray}\label{b613}
\begin{split}
C_{2100}=\frac{1}{6}q_{1}(0)F_{2100}\gamma_{40},~C_{1011}=\frac{1}{6}q_{1}(0)F_{1011}\gamma_{22},\\
C_{0021}=\frac{1}{6}q_{3}(0)F_{0021}\gamma_{04},~ C_{1110}=\frac{1}{6}q_{3}(0)F_{1110}\gamma_{22}.
\end{split}
\end{eqnarray}
where $\gamma_{ij}=\int_{0}^{l\pi}\gamma_{k_{1}}^{i}(x)\gamma_{k_{2}}^{j}(x)dx ~~~(i+j=4)$
and
\begin{eqnarray*}
&&F_{2100}=6\tau^{\ast}\big[3\tilde{\alpha}_{3}p_{12}^{2}\bar{p}_{12}+\tilde{\alpha}_{4}(p_{12}^{2}+2p_{12}\bar{p}_{12}),0\big]^{T},\\
&&F_{1011}=6\tau^{\ast}\big[6\tilde{\alpha}_{3}p_{12}p_{32}\bar{p}_{32}
+2\tilde{\alpha}_{4}(p_{12}p_{32}+p_{12}\bar{p}_{32}+p_{32}\bar{p}_{32}),0\big]^{T},\\
&&F_{0021}=6\tau^{\ast}\big[3\tilde{\alpha}_{3}p_{32}^{2}\bar{p}_{32}+\tilde{\alpha}_{4}(p_{32}^{2}+2p_{32}\bar{p}_{32}),0\big]^{T},\\
&&F_{1110}=6\tau^{\ast}\big[6\tilde{\alpha}_{3}p_{12}\bar{p}_{12}p_{32}
+2\tilde{\alpha}_{4}(p_{12}\bar{p}_{12}+p_{12}p_{32}+p_{32}\bar{p}_{12}),0\big]^{T}.
\end{eqnarray*}
Secondly, we have
\begin{eqnarray}\label{c613}
\begin{split}
D_{2100}&=\frac{1}{6l\pi}\Big(\frac{2}{-i\omega_{k_{1}}^{+}\tau^{\ast}}q_{1}^{2}(0)F_{2000}F_{1100}
+\frac{1}{i\omega_{k_{1}}^{+}\tau^{\ast}}q_{1}^{2}(0)F_{1100}F_{2000} \\
&+\frac{1}{i\omega_{k_{1}}^{+}\tau^{\ast}}q_{1}(0)\bar{q}_{1}(0)F_{1100}^{2}
+\frac{2}{3i\omega_{k_{1}}^{+}\tau^{\ast}}q_{1}(0)\bar{q}_{1}(0)F_{0200}F_{2000} \\
&-\frac{1}{i\omega_{k_{2}}^{-}\tau^{\ast}}q_{1}(0)q_{3}(0)F_{1010}F_{1100}
+\frac{1}{2i\omega_{k_{1}}^{+}\tau^{\ast}-i\omega_{k_{2}}^{-}\tau^{\ast}}q_{1}(0)q_{3}(0)F_{0110}F_{2000} \\
&+\frac{1}{i\omega_{k_{2}}^{-}\tau^{\ast}}q_{1}(0)\bar{q}_{3}(0)F_{1001}F_{1100}
+\frac{1}{2i\omega_{k_{1}}^{+}\tau^{\ast}+i\omega_{k_{2}}^{-}\tau^{\ast}}q_{1}(0)\bar{q}_{3}(0)F_{0101}F_{2000}
\Big),
\end{split}
\end{eqnarray}
where the coefficient vectors $F_{q_{1}q_{2}q_{3}q_{4}}(q_{1}+q_{2}+q_{3}+q_{4}=2)$ have the following forms.
\begin{eqnarray*}
\begin{array}{l}
F_{2000}=2\tau^{\ast}\left(\begin{array}{cc}
\tilde{\alpha}_{1}p_{12}-ae^{-i\omega_{k_{1}}^{+}\tau^{\ast}}
+\tilde{\alpha}_{2}p_{12}^{2} \\
c p_{12}-mp_{12}^{2}  \\
\end{array}\right),\\
F_{1100}=2\tau^{\ast}\left(\begin{array}{cc}
\tilde{\alpha}_{1}(p_{12}+\bar{p}_{12})-a(e^{i\omega_{k_{1}}^{+}\tau^{\ast}}
+e^{-i\omega_{k_{1}}^{+}\tau^{\ast}})+2\tilde{\alpha}_{2}p_{12}\bar{p}_{12} \\
c(p_{12}+\bar{p}_{12})-2mp_{12}\bar{p}_{12}  \\
\end{array}\right),
\end{array}
\end{eqnarray*}
$F_{1010}$, $F_{1001}$, $F_{0200}$, $F_{0110}$, $F_{0101}$, $F_{0020}$, $F_{0011}$ and $F_{0002}$ are omitted here. $D_{1011}$, $D_{0021}$, $D_{1110}$ could also be obtained, and more
details please refer to \citep{Du2018Double}.

Finally, we obtain
\begin{eqnarray}\label{603d}
\begin{split}
E_{2100}&=\frac{1}{6\sqrt{l\pi}}q_{1}(0)[S_{yz_{1}}(w_{01100})+S_{yz_{2}}(w_{02000})], \\
E_{1011}&=\frac{1}{6\sqrt{l\pi}}q_{1}(0)[S_{yz_{1}}(w_{00011})+S_{yz_{3}}(w_{01001})+S_{yz_{4}}(w_{01010})], \\
E_{0021}&=\frac{1}{6\sqrt{l\pi}}q_{3}(0)[S_{yz_{3}}(w_{00011})+S_{yz_{4}}(w_{00020})], \\
E_{1110}&=\frac{1}{6\sqrt{l\pi}}q_{3}(0)[S_{yz_{1}}(w_{00110})+S_{yz_{2}}(w_{01010})+S_{yz_{3}}(w_{01100})], \\
\end{split}
\end{eqnarray}
where $S_{yz_{i}} (i=1,2,3,4)$ are linear operators and
\begin{eqnarray}\label{g614}
S_{yz_{i}}(\varphi)=F_{y(0)z_{i}}\varphi(0)+F_{y(-1)z_{i}}\varphi(-1),
\end{eqnarray}
\begin{eqnarray*}
\begin{array}{l}
F_{y(0)z_{1}}=2\tau^{\ast}\left(\begin{array}{ll}
\tilde{\alpha}_{1}p_{12}-ae^{-i\omega_{k_{1}}^{+}\tau^{\ast}} & \tilde{\alpha}_{1}+2\tilde{\alpha}_{2}p_{12} \\
c p_{12} & c-2mp_{12} \\
\end{array}\right),
F_{y(-1)z_{1}}=2\tau^{\ast}\left(\begin{array}{ll}
-a & 0 \\
0 & 0 \\
\end{array}\right),\\
F_{y(0)z_{3}}=2\tau^{\ast}\left(\begin{array}{ll}
\tilde{\alpha}_{1}p_{32}-ae^{-i\omega_{k_{2}}^{-}\tau^{\ast}} & \tilde{\alpha}_{1}+2\tilde{\alpha}_{2}p_{32} \\
c p_{32} & c-2mp_{32} \\
\end{array}\right),\\
F_{y(-1)z_{1}}=F_{y(-1)z_{2}}=F_{y(-1)z_{3}}=F_{y(-1)z_{4}},~
F_{y(0)z_{2}}=\overline{F_{y(0)z_{1}}},
F_{y(0)z_{4}}=\overline{F_{y(0)z_{3}}}.
\end{array}
\end{eqnarray*}
Then we get
\begin{eqnarray*}
\begin{split}
w_{02000}(0)&=\frac{1}{\sqrt{l\pi}\tau^{\ast}}\Big[\frac{1}{-i\omega_{k_{1}}^{+}}p_{1}(0)q_{1}(0)
-\frac{1}{3i\omega_{k_{1}}^{+}}\bar{p}_{1}(0)\bar{q}_{1}(0)
+\frac{1}{(i\omega_{k_{2}}^{-}-2i\omega_{k_{1}}^{+})}p_{3}(0)q_{3}(0)\\&-
\frac{1}{(2i\omega_{k_{1}}^{+}+i\omega_{k_{2}}^{-})}\bar{p}_{3}(0)\bar{q}_{3}(0)
-
[-2i\omega_{k_{1}}^{+}I_{d}+\tilde{A}_{0}+
\tilde{B}_{0}e^{-2i\omega_{k_{1}}^{+}\tau^{\ast}}]^{-1}\Big]F_{2000},\\
w_{02000}(\theta)&=\frac{1}{\sqrt{l\pi}\tau^{\ast}}\Big[\frac{1}{-i\omega_{k_{1}}^{+}}p_{1}(\theta)q_{1}(0)
-\frac{1}{3i\omega_{k_{1}}^{+}}\bar{p}_{1}(\theta)\bar{q}_{1}(0)\\
&+\frac{1}{(i\omega_{k_{2}}^{-}-2i\omega_{k_{1}}^{+})}p_{3}(\theta)q_{3}(0)-
\frac{1}{(2i\omega_{k_{1}}^{+}+i\omega_{k_{2}}^{-})}\bar{p}_{3}(\theta)\bar{q}_{3}(0)\\
&-
e^{2i\omega_{k_{1}}^{+}\tau^{\ast}\theta}[-2i\omega_{k_{1}}^{+}I_{d}+\tilde{A}_{0}+
\tilde{B}_{0}e^{-2i\omega_{k_{1}}^{+}\tau^{\ast}}]^{-1} \Big]F_{2000}.
\end{split}
\end{eqnarray*}
Similarly, the expressions of $w_{01100}(\theta)$, $w_{00011}(\theta)$, $w_{01001}(\theta)$, $w_{01010}(\theta)$,
$w_{00020}(\theta)$ and $w_{00110}(\theta)$ could be given.

\end{document}